\documentclass[11pt,a4paper]{amsart} 
\usepackage{amssymb,xspace,amsfonts,enumerate,graphics}
\usepackage{amstext}
\theoremstyle{plain}

\date{\today}
\title{Extension of holomorphic functions defined on singular complex hypersurfaces with growth estimates in strictly pseudoconvex domains of $\cc^n$}
\author{William ALEXANDRE}
\address{Laboratoire Paul Painlev\'e U.M.R. CNRS 8524, U.F.R. de
Math\'ematiques,  cit\'e scientifique, Universit\'e Lille 1, F59 655 Villeneuve d'Ascq Cedex, France.}
\email{ william.alexandre@math.univ-lille1.fr}
\author{Emmanuel MAZZILLI}
\address{Laboratoire Paul Painlev\'e U.M.R. CNRS 8524, U.F.R. de
Math\'ematiques,  cit\'e scientifique, Universit\'e Lille 1, F59 655 Villeneuve d'Ascq Cedex, France.}
\email{ emmanuel.mazzilli@math.univ-lille1.fr}
\date{}
\newtheorem{theorem}{Theorem}[section]
\newtheorem{lemma}[theorem]{Lemma}

\newtheorem{proposition}[theorem]{Proposition}
\newtheorem{corollary}[theorem]{Corollary}
\newtheorem{remark}{\it Remark}
\newtheorem{definition}[theorem]{Definition}
\def \pint {\vbox{ \hbox to 5 pt {\hfil \vrule height 4pt}\hrule}\hskip 3pt}
\def\tq{/\ }
\def\leqs{\lesssim}
\def\geqs{\gtrsim}
\def\eqs{\eqsim}
\def\cc{\mathbb{C}}
\def\rr{\mathbb{R}}
\def\nn{\mathbb{N}}
\def\zz{\mathbb{Z}}
\def\re{{\rm Re} }
\def\im{{\rm Im} }
\def \qed {\hbox{\hskip 5pt} \vbox{\hrule \hbox to 5pt 
{\vrule height 4.2pt \hfil \vrule}\hrule}}
\def \pint {\vbox{ \hbox to 5 pt {\hfil \vrule height 4pt}\hrule}\hskip 3pt}

\newcommand{\cal}{\mathcal}

\newcommand{\diffp}[2]{\frac{\partial #1}{\partial #2}}
\newcommand{\mlabel}[1]{\label {#1}}
\renewcommand{\over}[2]{\genfrac{}{}{0pt}{}{#1}{#2}}

\newcommand{\p}[2]{{\cal P}_{#1}(#2)}

\newcommand{\pr} {\noindent{\it Proof:} }
\newcommand{\ko} {Koranyi }
\newcommand{\kb}[2] {{\cal P}_{#1}(#2)}
\newcommand{\kbz}[1] {{\cal P}_{c|\rho(#1)|}(#1)}
\newcommand{\kbc}{\left(\kb{c|\rho(z_j)|}{z_j}\right)_{j\in\nn}}
\newcommand{\cu}{{\cal U}}
\newcommand{\cv}{{\cal V}}
\newcommand{\co}{{\cal O}}
\def\sing{{\rm Sing}}
\def\reg{{\rm Reg}}
\def\dim{{\rm dim\, }}
\begin{document}
\pagestyle{plain}

\begin{abstract}
Let $D$ be a strictly pseudoconvex domain and $X$ be a singular analytic set of pure dimension $n-1$ in $\cc^n$ such that $X\cap D\neq \emptyset$ and $X\cap bD$ is transverse. We give sufficient conditions for a function holomorphic on $D\cap X$ to admit  a holomorphic extension which belongs to $L^q(D),$ $q\in [1,+\infty[$, or to $BMO(D)$. The extension is given by mean of integral representation formulas and residue currents.
\end{abstract}
\maketitle
\section{Introduction}
In this article we are interested in the following natural question : Given an analytic set $X$, a strictly pseudoconvex domain $D$ such that $X\cap D\neq\emptyset$ and a function $h$ holomorphic on $X\cap D$, does their exist a function $H$, holomorphic on $D$, such that the restriction of $H$ to $X\cap D$ coincides with $h$ ? Cartan's theorem B asserts that the answer is positive, even if $D$ is only pseudoconvex. Difficulties arise when we want $H$ to satisfy growth conditions like being bounded or belonging to $L^q(D)$ or $BMO(D)$. 

 When $X$ is a hyperplane and $D$ is pseudoconvex, Ohsawa and Takegoshi  proved in \cite{OT} that any $h\in L^2(X\cap D)\cap {\cal O}(X\cap D)$ admits an extension $H\in L^2(D)\cap \co (D)$. Ohsawa in \cite{OhsII} generalized this result to complex manifolds of higher codimension. When $D$ is strictly pseudoconvex, Amar in \cite{Ama} proved that any function holomorphic and bounded on $X\cap D$ has an extension holomorphic and bounded on $D$. The case of bounded functions defined on a manifold and extended to bounded functions on weakly pseudoconvex was also positively solved in the case of convex domain of finite type in \cite{Ale1} and \cite{DiM2}.

When $X$ is no longer a manifold but a singular variety, things are much more complicated. In \cite{DiM0}, Diederich and the second author exhibit, in $\cc^3$, an algebraic complex hypersurface $X$ with singularities and a function $h$ bounded and holomorphic in the intersection  of $X$ and the unit ball which does not have a $L^2$-extension in the unit ball.  This result was generalized to weighted $L^2$-spaces by Guan and Li in \cite{GL}.

In \cite{Ber}, for $X$ a singular varieties and $D$ pseudoconvex, Berndtsson gave a condition under which $h\in\co(X\cap D)$ has an extension in $\co(D)\cap L^2(D)$. However if $h$ satisfies this condition, $h$ must vanish on $\sing(X)$, the set of singular points of $X$. Hence, if $h$ is constant on $X\cap D$, Berndtsson's result does not provide a holomorphic extension of $h$ which belongs to $L^2(D)$ while it trivially exists. In \cite{AM}, when $X$ is a singular variety in $\cc^n$, when $D$ is a strictly convex domain and when $q$ belongs to $[1,+\infty],$ the authors gave necessary conditions for $h\in \co(X\cap D)$ to have an extension $H\in \co(D)\cap L^q(D)$. When $n=2$ and $q\in [1,+\infty[$, they also proved that these conditions are sufficient and when $q=+\infty$, they proved that these conditions imply that $H$ belongs to $\co (D)\cap BMO(D)$. In this article, we aim at generalizing the results of \cite{AM} for $n\geq 2$ and for strictly pseudoconvex domains. Moreover, we want to give conditions which are stable by biholomorphic transformations.
\par\bigskip

Let $D$ be a strictly pseudoconvex domain in $\cc^n$ with smooth boundary. We write $D$ as $D=\{z\in\cu \tq \rho(z)<0\}$ where $\rho$ is a smooth strictly plurisubharmonic function  defined in a neighborhood $\cu$ of $\overline D$, such that the gradient of $\rho$ does not vanish in a neighborhood of $bD$, the boundary of $D$. We denote by $D_t$, $t\in\rr$, the set $D_t=\{z\in\cc^n,\ \rho(z)<t\}$, by $bD_t$ its boundary, by $\eta_p$  the outer 
unit normal to $bD_{\rho(p)}$ at a point $p\in\cu$ and by $T_p^\cc bD_{\rho(p)}$ the complex tangent space to $bD_{\rho(p)}$ at the point $p$. For all $p$ in $\cu$, we denote by $\epsilon_1(p)=\eta_p,\epsilon_2(p),\ldots,\epsilon_n(p)$ an orthonormal basis of $T^\cc_pbD_{\rho(p)}$ which locally smoothly depends on $p$.

\par\smallskip

Let $X=\{\zeta\in\cu \tq f(\zeta)=0\}$ be an analytic set of pure dimension $n-1$ in $\cc^n$. 
We  denote by $\reg(X)$ the set of regular points of $X$ and by $\sing(X)$ the set of singular points of $X$.
We also denote by $C_4(X,p)$ the fourth Whitney tangent cone at $p\in X$ to $X$, that is the set of vectors $v$ for which there are sequences of points $(z_j)_j$ of $\reg(X)$ converging to $p$ and vectors $(v_j)_j$ converging to $v$ such that, for all $j$,  $v_j$ is tangent to $X$ at $z_j$. For all $p\in X$, $\dim C_4(X,p)\geq n-1$ and Stutz showed in \cite{Stu} that the algebraic set $J=\{p\in\sing(X)/\ \dim C_4(X,p)> n-1\}$ has dimension at most $n-2$ and is thus strictly smaller than $\sing(X)$. In this paper, we will assume that $J$ is in fact empty. 


\par\smallskip

In all this work, we assume that $X$ and $D$ satisfy the following assumptions :
\begin{enumerate}[(a)]
 \item \label{hyp1} The intersection $X\cap D$ is non empty.
 \item \label{hyp2} The intersection $X\cap bD$ is transverse in the sense of tangent cones, that is, for all $p\in bD\cap X$, the vector space generated by $T_p^\cc bD$ and $C_4(X,p)$ is $\cc^n$.
 \item \label{hyp3} $X$ is irreducible at all point $p\in X\cap bD$.
 \item  \label{hyp4} Any point $p\in bD\cap X$ is either a regular point of $X$ or $p$ is a regular point of $\sing(X)$, such that  $\dim C_4(X,p)=n-1$.
\end{enumerate}

We will prove the following theorem which generalizes Theorem 1.1 of \cite{AM} :
\begin{theorem}\label{th0}
Let $X$ be an analytic set of pure dimension $n-1$ in $\cc^n$ and let $D$ be a strictly pseudoconvex domain in $\cc^n$  which satisfy (\ref{hyp1}-\ref{hyp4}).\\
There exists an integer $k\geq 1$ depending only on $X$ such that if $h$ is a holomorphic function on $X\cap D$ which has a $C^\infty$ smooth extension $\tilde h$ on $D$ which satisfies
\begin{enumerate}[(i)]
 \item \label{th0i} there exists $N\in\nn$ such that $|\rho|^N \tilde h$ vanishes at order $k$ on $bD$,
 \item \label{th0ii}there exists $q\in[1,+\infty]$ such that  $\left|\diffp{^{|\alpha|}\tilde h}{\overline{\epsilon_1}^{\alpha_1}\ldots \partial \overline{\epsilon_n}^{\alpha_n}}\right||\rho|^{\alpha_1+\frac{\alpha_2+\ldots+\alpha_n}2}$ belongs to $L^q(D)$ for all multi-index $\alpha$  
 with $|\alpha|\leq k$,
 \item\label{th0iii} $\diffp{^{|\alpha|}\tilde h}{\overline{\epsilon_1}^{\alpha_1}\ldots \partial \overline{\epsilon_n}^{\alpha_n}}=0$ on $X\cap D$ for all
 multi-index $\alpha$ with $0<|\alpha|\leq k$,
\end{enumerate}
then $h$ has a holomorphic extension $H$ in $L^q(D)$ when $q<+\infty$ and in $BMO(D)$ when $q=+\infty$. Moreover, up to a uniform multiplicative constant
depending only on $k$ and $N$, the norm of $H$ is  bounded 
by the supremum of the $L^q$-norm of $\zeta\mapsto
\left|\diffp{^{\alpha}\tilde h}{\overline{\epsilon_1}^{\alpha_1}\ldots \partial \overline{\epsilon_1}^{\alpha_n}}\right||\rho|^{\alpha_1+\frac{\alpha_2+\ldots+\alpha_n}2}$ for $\alpha$ multi-index such that $|\alpha|\leq k$.
\end{theorem}
The extension will be given as in \cite{AM} by an integral operator combining a Berndtsson-Andersson operator and the $\overline \partial$ of a current $T$ which is a refinement of the current  of \cite{Maz1} and which satisfies $fT=1$. The current $\overline\partial T$ is some kind of a perturbation of  the classical current $\overline\partial \left[\frac1f\right]$ defined in  \cite{HL}. The former current is defined using Hironaka's theorem and thus is not very explicit. In order to have precise estimates of the extension, we need a completely explicit current. We will defined it in a neighborhood of any point $p\in D$ by taking into account of branches of $X$ which are ``close'' or ``far'' from $p$.
Having tools which enable us to quantify how a branch is close or far from a point, is the first problem we have to overcome.
The first part of the solution will be given by Koranyi balls and the structure of homogeneous space of a strictly pseudoconvex domains which will provide a good (pseudo-)metric in order to quantify the meaning of being far or close from a point $p$. The second part of the answer will be given by Hypothesis (\ref{hyp3}-\ref{hyp4}) which provide us a parametrization of $X$. 

The current will act by differentiation on a Berndtsson-Andersson kernel and such kernels have a worse behavior when they are differentiated in the normal direction than in the tangential directions to the boundary of the domain. In dimension $2$, there is only one direction but not when $n>2$ and we have to find the good one. Moreover, in dimension 2,  the singularities of an analytic set are always isolated. In particular, a singularity is either far from the boundary of the domain or belongs to the boundary of the domain and so, there is a kind of dichotomy. When $n>2$, there can be a ``continuum'' of singularity which goes through the boundary of $D$.  It should be noticed that in most papers dealing with singular varieties and giving fine estimates, singularities are often assumed isolated, but we do not make such an assumption here and we will learn to deal with non isolated singularities. 

We will face an other difficulty which does not appear in \cite{AM} and which is due to the non convexity of $D$.  In order to define the extension operator, we have to write $f(z)-f(\zeta)=\sum_{j=1}^nb_j(\zeta,z)(z_j-\zeta_j)$. Any such decomposition $b=\sum_{j=1}^nb_jd\zeta_j$ enables us to define an extension operator but not all $b$ give a good extension because we need some kind of compensation between $b$ and the term $\frac1f$  which will appear in the kernel of our operator (see Lemma \ref{lemme3.1} for details). In \cite{AM}, the convexity of $D$ was used in order to have $b_j$ equal to $\sum_{|\alpha|\leq k} \frac1{(|\alpha|+1) \alpha!}\diffp{^{|\alpha|+1}f }{\zeta_j \partial \zeta^\alpha}(\zeta)(z-\zeta)^\alpha+O(|\zeta-z|^{k+1})$. Here, we will prove that given any Hefer decomposition of $f$, we can construct a good decomposition which will be equal to the derivatives of $f$ up to order $k$ (see Lemma \ref{lemme Hefer}).
\par\medskip
In Theorem \ref{th0}, we assume the existence of a smooth extension $H$ of $h\in\co(X\cap D)$, although no such assumption was made in the previous papers (see \cite{Ama,Cum} for example). However, during the proof of the existence of good extensions, often a smooth extension is first constructed. In general, this construction is a bit obvious : if $h$ is defined on a manyfold, maybe after a local biholomorphism, this manyfold is equal to the set $\{\zeta_1=0\}$ and the smooth extension is just $\tilde h(\zeta)= h(0,\zeta_2,\ldots,\zeta_n)$. This smooth extension will then satisfies the hypothesis (\ref{th0i}-\ref{th0iii}) of Theorem \ref{th0}. Of course, a singular variety cannot be written in such a way and, in our case, we have to assume the existence of such a smooth extension. The next question is thus ``when such a smooth extension does exist ?''.

In \cite{AM} we gave conditions under which such extensions exist in $\cc^2$. They were formulated in terms of control of divided differences on linear discs. Roughly speaking, a divided difference of order $m$ of a holomorphic function $H$ near a singularity is close to a $m$-th derivative of $H$ and Cauchy inequalities tell us that if $H$ is bounded in $D$ or if $H$ belongs to $L^q(D)$, we control its derivatives and thus its divided differences with powers of $d(\cdot,bD)$, the distance to the boundary of $D$. These conditions, even if necessary and sufficient or nearly sufficient, were not stable under biholomorphism. We want to give here conditions which are kept under biholomorphic transformations. Let $\Delta$ be the unit disc of $\cc$. We will prove when $q=+\infty$ the following theorem. 

\begin{theorem}\label{th1}
Let $X$ be an analytic set of pure dimension $n-1$ in $\cc^n$ and let $D$ be a strictly pseudoconvex domain in $\cc^n$  which satisfy (\ref{hyp1}-\ref{hyp4}). Let $h$ belongs to $\co(X\cap D)$. If for all holomorphic disc $\gamma : \Delta\to D$ such that $\gamma(\Delta)\cap X\neq \emptyset$, there exists $h_\gamma\in\co(\Delta)$ such that
\begin{itemize}
 \item $\sup_\Delta |h_\gamma|$ is bounded uniformly with respect to $\gamma$,
 \item for all $t\in \Delta$ such that $\gamma(t)$ belongs to $X$, $h_\gamma(t)=h\circ\gamma(t)$,
\end{itemize}
then there exists $\tilde h$ which satisfies the assumptions (\ref{th0i}-\ref{th0iii}) of Theorem \ref{th0}.
\end{theorem}

Of course, if $h$ has a bounded holomorphic extension, it satisfies the assumptions of Theorem \ref{th1} and these assumptions are kept under biholomorphic transformations. 

Since any $h$ belonging to $\co(X\cap D)$ always has a holomorphic extension to $D$, it suffices in Theorem \ref{th1} to consider not all discs but  holomorphic discs close to the boundary of $D$. We will also see that it suffices to consider only regular holomorphic discs $\gamma$ such that the intersection of $\gamma(\Delta)$ and any branch of $X$ is either a singleton or empty. And in fact, provided we are in suitable coordinates, it suffices to consider linear discs, i.e.\! images of $\Delta$ by complex affine transformations. 

Theorems \ref{th0} and \ref{th1} together assert that if $h\in\co(X\cap D)$ can be extended holomorphically to  holomorphic discs with a uniform bound, then $h$ has a holomorphic extension to $D$ which belongs to the class $BMO(D)$. When $X$ is a manyfold, and when $h$ is holomorphic and bounded on $X\cap D$, given a holomorphic disc $\gamma$ such that $\gamma(\Delta)\cap X$ is a singleton, we can trivially extend $h$ to $\gamma(\Delta)$ and we thus get a uniformly bounded extension of $h$ to $\gamma(\Delta)$. We can therefor apply Theorems \ref{th1} and \ref{th0} and we get a holomorphic extension of $h$ to $D$ which belongs to the class $BMO(D)$ and we thus nearly recover Henkin's results in \cite{Hen}.

In Theorem \ref{th2}, still using extension on holomorphic discs, we will also give an analog of Theorem \ref{th1} for $L^q$-extension when $q<+\infty$. However in order to have $L^q$-extension, we need some kind of average properties that will be formulated with homogeneous covering and Koranyi discs. The condition of Theorem \ref{th2} will also be stable under biholomorphic transformations and will be necessary and sufficient for $h$ to have a holomorphic extension which belongs to $L^q(D)$.

The paper is organized as follows. In  Section \ref{section1}, we construct the extension operator. In Section \ref{section2}, we prove that it satisfies the conclusion of Theorem \ref{th0}. In Section \ref{section3}, we prove Theorem \ref{th1} and give necessary and sufficient conditions on $h$ for $h$ to have a holomorphic extension which belongs to $L^q(D)$, $q<+\infty$.

\section{Construction of the current}\label{section1}
We first want to define a current $T$ such that $fT=1$. We will define $T$ locally but we first need more analytic informations on the analytic set $X$. We denote by $d(z,bD)$ the distance from $z$ to $bD$ and we define the anisotropic  Koranyi balls centered at $z$ of radius $r>0$ by $\kb rz:=\{z+\lambda \eta_z+\mu v\tq v\in T^\cc_zbD_{\rho(z)},\ |\lambda|+ |\mu|^2< r\}.$\\
We also adopt the following convention : we will write $A\leqs B$ if there exists a constant $c>0$ such that $A\leq cB$. We write $A\eqs B$ if both $A\leqs B$ and $B\leqs A$ hold true. Moreover, through out this paper, $C$ will always represent a big constant and $c$ a small one such that choosing $c$ smaller does not imply that $C$ must be chosen bigger. However, we accept that choosing $C$ bigger implies that $c$ must be chosen smaller.

\subsection{Local parametrization of $X$}\label{Local parametrization}
Without restriction, we assume that $0$ belongs to $\sing(X)\cap bD$ and we work near 0.\\
Hypothesis (\ref{hyp3}) and (\ref{hyp4}) and Proposition 4.2 of \cite{Stu} imply that 
there are a neighborhood $\cu(0)$ and a one to one holomorphic map $\Phi$ of a neighborhood $\cv(0) \subset \cc^{n-1}$ of $0$ onto $\cu(0)\cap X$ such that 
 \begin{enumerate}[(i)]
  \item $\sing(X)\cap \cu(0)=\Phi(\cv(0)\cap\{\zeta_1=0\})$,
  \item $\Phi:\cv(0)\setminus \{t_1=0\}\to\reg(X)\cap \cu(0)$ is biholomorphic,
  \item \label{c} after perhaps a holomorphic change of coordinates, $\Phi$ is of the form
  $$\Phi(t)=(t^k_1,t_2,\ldots,t_{n-1},t^k_1\varphi(t)),$$
  where $k$ is the multiplicity of the cover of $X$ over $\cc^{n-1}$ in a neighborhood of $0$, $\varphi$ is holomorphic in $\cv(0)$ and $\varphi(0)=0$.
 \end{enumerate}
Thus, locally and maybe after a local change of coordinates, $\sing(X)=\{z\in\cc^n\tq z_1=z_n=0\}$ and by transversality, we deduce that $\eta_0\neq (0,\ldots, 0,1)$.

Moreover we can write $X$ near $0$ as 
$$X=\left\{z\in\cc^n\tq P(z)=\prod_{j=0}^{k-1}\left(z_n-z_1\varphi(z_1^{\frac1k} \omega^j,z_2,\ldots,z_{n-1})\right)\right\}$$
where $z_1^{\frac1k}$ satisfies $(z_1^{\frac1k})^k=z_1$. 
\par\medskip
Intuitively, in these local coordinates, the action of the current $T$ on a test function $\phi$ will be of the kind $\langle T,\phi\rangle=\int \frac{\overline{P}(\zeta)}{f(\zeta)}\diffp{^k\phi}{\overline \zeta_n^k}$. Since $f=uP$ with $u$ zero free, we integrate a bounded quantity and so $T$ is well defined. When we integrate $k$ times by parts, we get $fT=1$.

The extension operator will make $T$ acting on a Berndtsson-Andersson reproducing kernel. Thus we will differentiate  a reproducing kernel $k$ times and it is well known that, near the boundary of a domain, these derivatives explode like a power of $\frac1{d(z,bD)}$. More precisely, a derivative in the normal direction implies a loss of a factor $d(\cdot,bD)$ and a derivative in a tangent direction implies a loss of a factor  $d(\cdot,bD)^{\frac12}$. Thus, our interest will be to differentiate only in a tangential direction. In order to do this, we will cover a neighborhood of $bD$ with Koranyi balls $\p {d(z_j,bD)} {z_j}$, $j\in\nn$, and on each Koranyi ball, we will define a current $T_j$ of the previous kind such that we differentiate in a tangential direction at $z_j$. Thus for any point $z$ close enough to 0, we will have to choose a basis depending on $z$, not necessarily orthonormal, such that for example the last vector of this basis is tangent to $bD_{\rho(z)}$ at $z$ and such that in the coordinates induces by this new basis, $X$ has a parametrization of type $t\mapsto (t_1^k,t_2,\ldots,t_{n-1},t_1^k\varphi_z(t))$.

If we denote by $(e_1,\ldots,e_n)$ the canonical orthonormal basis, we have $\eta_0\neq \pm e_n$. So, for $z$ near $0$, we can put $w_z=\frac1{|e_n-\langle e_n,\eta_z\rangle\eta_z|}(e_n-\langle e_n,\eta_z\rangle\eta_z)=w_{z,1}e_1+\ldots+w_{z,n} e_n$. Since $e_n\neq \eta_0$, we have $w_{0,n}\neq 0$ and, by continuity, $w_{z,n}\neq 0$ for all $z$ sufficiently close to $0$.  For such a point $z$, we denote by $A_z=[e_1,\ldots, e_{n-1},w_z]$ the change of coordinates matrix from the canonical basis to $(e_1,\ldots, e_{n-1},w_z)$ and for $\zeta\in\cc^n$, we put $\pi_z(\zeta)=(\pi_{z,1}(\zeta),\ldots,\pi_{z,n}(\zeta))=A^{-1}_z\zeta$. We want to prove the following lemma which gives us a parametrization of $X$ in the coordinates induced by the basis $(e_1,\ldots, e_{n-1},w_z)$.
 \begin{lemma}\label{lemme2.1}
  There exist a neighborhood  $\cu(0)$ of the origin in $\cc^{n}$ and a neighborhood $\cv(0)$ of the origin in $\cc^{n-1}$ such that for all $z$ close enough to $0$, there exists a function $\varphi_z$ holomorphic in $\cv(0)$ which satisfies
  \begin{enumerate}[(i)]
   \item $|\varphi_z(t')| \leqs |t'|$, uniformly with respect to $z$ and $t'$,
   \item $X\cap \cu(0)= \pi_z^{-1}\left\{\left({t'_1}^k,t'_2,\ldots,t'_{n-1},{t'_1}^k\varphi_z(t')\right)\tq t'=(t'_1,\ldots, t'_{n-1})\in\cv(0)\right\}$.
  \end{enumerate}
  Moreover, $\varphi_z$ can be assumed  uniformly bounded.
 \end{lemma}
\pr
 We fix a holomorphic $k$-th roots in the disc $D(1,1)\subset\cc$ and put $$\tilde\Phi:(z,t,t')\mapsto\left(\begin{array}{c}t'_1-t_1\left(1-\frac{w_{z,1}}{w_{z,n}} \varphi(t)\right)^{\frac1k}\\
				    t'_2-\left(t_2-\frac{w_{z,2}}{w_{z,n}} t_1^k \varphi(t)\right)\\
                                    \vdots\\
                                    t'_{n-1}-\left(t_{n-1}-\frac{w_{z,n-1}}{w_{z,n}} t_1^k\varphi(t)\right)\\
                                    \end{array}\right)$$
where $\varphi$ is given by (\ref{c}). Since $\varphi(0)=0$, $\diffp{\tilde \Phi}{t}(0,0,0)=-Id_{\cc^{n-1}}$. The implicit functions theorem implies that there exists a neighborhood $\cu(0)$ of $0\in\cc^n$, two neighborhoods $\cv(0)$ and $\cv'(0)$ of $0\in\cc^{n-1}$ and $\tilde \Psi:\cu(0)\times\cv'(0)\to\cv(0)$, holomorphic with respect to $t'$ and smooth with respect to $z$, such that for all $(z,t,t')\in\cu(0)\times\cv(0)\times\cv'(0)$, $\tilde\Phi(z,t,t')=0$ if and only if $t=\tilde\Psi(z,t')$.Moreover $\tilde\Psi(z,0)=0$ so $|\tilde\Psi(z,t')|\leqs|t'|$, uniformly with respect to $z$ and $t'$.\\
Differentiating the following equality
\begin{align}
t'_1-\tilde\Psi_1(z,t')\left(1-\frac{w_{z,1}}{w_{z,n}} \varphi\circ\tilde\Psi(z,t')\right)^{\frac1k}&=0, \label{eq1}
 \end{align}
 we get for all multi-index $\alpha=(0,\alpha_2,\ldots,\alpha_{n-1})$
 $$\diffp{^{|\alpha|}\tilde\Psi_1}{ {t'}^\alpha} (z,0)=0.$$
 Differentiating (\ref{eq1}) with respect to $t'_1$ and evaluating at $(z,0)$, we get
 $$\diffp{\tilde\Psi_1}{t'_1}(z,0)=1.$$
 Therefore there exists a holomorphic function $\tilde \varphi_z$ such that  $\tilde\Psi_1(z,t')=t'_1(1+\tilde\varphi_z(t'))$
 and $|\tilde \varphi_z(t')|\leqs |t'|$, uniformly with respect to $z$ and $t'$.\\
 Finally, we put $\varphi_z(t')=\frac1{w_{z,n}}\left(1+\tilde\varphi_z(t')\right)^k\varphi(\tilde\Psi(z,t'))$.
 
 If $\zeta$ is close to $0$ and belongs to $X$, then it can be written as $\zeta=(t_1^k,t_2,\ldots,t_{n-1},t_1^k\varphi(t))$. 
 Putting $t'=\left(t_1\left(1-\frac{w_{z,1}}{w_{z,n}} \varphi(t)\right)^{\frac1k}, t_2-\frac{w_{z,2}}{w_{z,n}} t_1^k \varphi(t),\ldots, t_{n-1}-\frac{w_{z,n-1}}{w_{z,n}} t_1^k\varphi(t)\right)$ we thus obtain  $\pi_z(\zeta)=(t'^k_1, t'_2,\ldots,t'_{n-1},t'^k_1 \varphi_z(t'))$ and conversely.\qed
 \par\bigskip
 For a given point $z$, we will denote by $\zeta'$ the coordinates of a point $\zeta$ in the coordinates system centered at 0 and of basis $e_1,\ldots,e_{n-1}, w_z$.  We define
 $$P_z(\zeta')=\prod_{j=0}^{k-1}\left(\zeta'_n-\zeta'_1\varphi_z(\omega^j{\zeta'_1}^{\frac1k},\zeta'_2,\ldots,\zeta'_{n-1})\right)$$
 where $\omega=e^{\frac{2i\pi}k}$ and ${\zeta'_1}^{\frac1k}$ is any complex number such that $({\zeta'_1}^{\frac1k})^k=\zeta'_1$.
 Therefor $\zeta$ belongs to $X$ if and only if $P_z(\pi_z(\zeta))=0$. We want to link uniformly $P_z$ to $f$ :
 \begin{proposition}\label{prop2.2}
 For all $z$ near $0$, there exists a holomorphic function $u_z$ such that $f(\zeta)=u_z(\zeta)P_z(\pi_z(\zeta))$ for all $\zeta$ in a neighborhood of $0$ which does not depends on $z$ and $|u_z|\eqs 1$ uniformly with respect to $z$.
 \end{proposition}
\pr   We consider the family of functions $(f_z)_z$ where for all $z$, $f_z(\zeta')=f(\pi_z^{-1}(\zeta'))=f(\zeta'_1e_1+\ldots+\zeta'_{n-1} e_{n-1}+\zeta'_n w_z)$. We first want to apply Rouch\'e's theorem to  $f_0(0,\ldots,0,\cdot)$ and $f_z(\zeta'_1,\ldots,\zeta'_{n-1},\cdot)$.\\
Weierstrass Preparation Theorem implies that, near $0$,  there exist a non vanishing holomorphic function $u$ and a Weierstrass polynomial $P(\zeta)=\zeta_n^k+\zeta_{n}^{k-1}a_1(\zeta_1,\ldots,\zeta_{n-1})+\ldots+a_{k}(\zeta_1,\ldots,\zeta_{n-1})$ such that $f=uP$. Since the function $\zeta_n\mapsto f(0,\ldots,0,\zeta_n)$ vanishes at order $k$ at $0$, we also have  $a_1(0,\ldots, 0)=\ldots=a_k(0,\ldots,0)=0$.
Moreover, we can write 
$$P(\zeta)=\prod_{j=0}^{k-1}\left(\zeta_n-\zeta_1\varphi(\zeta_1^{\frac1k} \omega^j,\zeta_2,\ldots,\zeta_{n-1})\right)$$
from which we get
\begin{align*}
 f_0(0,\ldots,0,t)
&=u(tw_{0})\prod_{j=0}^{k-1}\left(tw_{0,n}-w_{0,1}t\varphi((tw_{0,1})^{\frac1k} \omega^j,tw_{0,2},\ldots,tw_{0,n-1})\right)\\
&={t}^k \tilde u(t)
\end{align*}
where $|\tilde u(t)|=|u(t w_0)|\prod_{j=0}^{k-1}\left|w_{0,n}-w_{0,1}\varphi((tw_{0,1})^{\frac1k} \omega^j,tw_{0,2},\ldots,tw_{0,n-1})\right|\eqs1$ if $|t|$ is small enough.\\
Therefore there exists a small $\varepsilon >0$ such that for $t$ in  $\overline{D(0,\varepsilon)}=\{\xi\in\cc\tq |\xi|\leq\varepsilon\}$, $f_0(0,\ldots,0,t)=0$ if and only if $t=0$. In particular, $f_0(0,\ldots,0,\cdot)$ does not vanish on $bD(0,\varepsilon)$, the boundary of the disc $D(0,\varepsilon)$. We put $a=\inf_{bD(0,\varepsilon)} |f_0(0,\ldots,0,\cdot)|>0$.\\
The mean value inequality then gives us $\varepsilon'>0$ such that if $(t_1,\ldots,t_{n-1},t)$ belongs to $D(0,\varepsilon')^{n-1}\times D(0,\varepsilon)$ then $|f_0(t_1,\ldots,t_{n-1},t)-f_0(0,\ldots,0,t)|<\frac a4$.

Again the mean value inequality gives us $\varepsilon''>0$ such that if $z$ belongs to $B(0,\varepsilon'')=\{z\in\cc^n\tq |z|<\varepsilon''\}$ and $(t_1,\ldots,t_{n-1},t)$ belongs to $D(0,\varepsilon')^{n-1}\times D(0,\varepsilon)$, then $|f_z(t_1,\ldots,t_{n-1},t)-f_0(t_1,\ldots,t_{n-1},t)|<\frac a4$. This yields $|f_z(t_1,\ldots,t_{n-1},t)-f_0(0,\ldots,0,t)|<\frac a2$
 for all $z\in B(0,\varepsilon'')$ and $(t_1,\ldots,t_{n-1},t)\in D(0,\varepsilon')\times D(0,\varepsilon)$.
 Rouch\'e's theorem then implies that $f_z(t_1,\ldots,t_{n-1},\cdot)$ has exactly $k$ zeros in the disc $D(0,\varepsilon)$.
 
 Now we apply Weierstrass Preparation Theorem to $f_z$ and we write $f_z=u_zQ_z$ where $Q_z(\zeta')={\zeta'_n}^k+a_1^{(z)}({\zeta'_1},\ldots,\zeta'_{n-1}){\zeta'_n}^{k-1}+\ldots +a_k^{(z)}(\zeta'_{1},\ldots,\zeta'_{n-1})$ is a Weierstrass polynomial and $u_z$ does not vanish in $D(0,\varepsilon')^{n-1}\times D(0,\varepsilon)$.
 Since $f_z({t'_1}^k,t'_2,\ldots,t'_n,{t'_1}^k \varphi_z(t'))=0$, we conclude that $Q_z=P_z$.\\
 For all $\zeta'\in D(0,\varepsilon')^{n-1}\times D(0,\varepsilon)$, we have
 \begin{align*}
  \left|\frac1{u_z(\zeta')}\right|&\leq \sup_{|t|=\varepsilon} \frac{|Q_z(\zeta'_1,\ldots,\zeta'_{n-1},t)|}{|f_z(\zeta'_1,\ldots,\zeta'_{n-1},t)|}.
 \end{align*}
 On the first hand, for all $(\zeta'_1,\ldots,\zeta'_{n-1})\in D(0,\varepsilon')^{n-1}$ and all $t\in bD(0,\varepsilon)$, we have
 \begin{align*}
  |f_z(\zeta'_1,\ldots,\zeta'_{n-1},t)|&\geq |f_0(0,\ldots, 0,t)|-|f_z(\zeta'_1,\ldots,\zeta'_{n-1},t)-f_0(0,\ldots, 0,t)|\geq \frac a2.
 \end{align*}
On the other hand, since $(f_z)_z$ converges uniformly to $f_0$ when $z$ tends to $0$, Weierstrass Preparation Theorem also implies that $(Q_z)_z$ converges uniformly to $Q_0$ when $z$ tends to $0$ and so $Q_z$ is uniformly bounded. Therefore we have $|u_z({\zeta'})|\geq 1$ uniformly with respect to $z$ and ${\zeta'}$ . Analogously, we also have $|u_z|\leq 1$ and so $|u_z|\eqs 1$.\qed
\par\medskip

We now want to understand the interplay of the geometries of $X$ and  $D$. This will be the goal of the following propositions, firstly near $\sing(X)$.

\begin{proposition}\mlabel{prop2.3}
There exists $C>0$ big enough such that for all $c>0$, all $z$ sufficiently close to $0$ such that $d(z,\sing(X))\leq 10C(c|\rho(z)|)^{\frac12}$, all $\zeta\in \kbz z$ and all $\xi\in\cc$ such that $\xi^k=\pi_{z,n}(\zeta)$, 
we have $\left|\pi_{z,n}(\zeta)-\pi_{z,1}(\zeta)\varphi_z(\xi,\pi_{z,2}(\zeta),\ldots,\pi_{z,n-1}(\zeta)\right|\leq C^2|c\rho(z)|^{\frac12}$ and 
the point $\pi^{-1}_z\big(\pi_{z,1}(\zeta),\ldots,\pi_{z,n-1}(\zeta), \pi_{z,1}(\zeta)\varphi_z(\xi,\pi_{z,2}(\zeta),\ldots,\pi_{z,n-1}(\zeta))\big)=\zeta+\big(\pi_{z,1}(\zeta)\varphi_z(\xi,\pi_{z,2}(\zeta),\ldots,\pi_{z,n-1}(\zeta))-\pi_{z,n}(\zeta)\big)w_z$ belongs to $X\cap \kb {cC^4|\rho(z)|}z$.
\end{proposition}
\pr If $d(z,\sing(X))\leq 10C|c\rho(z)|^{\frac12}$ then for all $\zeta\in \kbz z$ we have 
\begin{align*}
 d(\zeta,\sing(X))&\leq |\zeta-z|+d(z,\sing(X))\\
 &\leq (10C+1)|c\rho(z)|^{\frac12}.
\end{align*}
Since $d(\zeta ,\sing(X))=\sqrt{|\zeta _1|^2+|\zeta _n|^2}$, it comes
$ |\pi_{z,1}(\zeta)|=\left|\zeta _1-\frac{w_{z,1}}{w_{z,n}}\zeta _n\right|\leqs C|c\rho(z)|^{\frac12},$ and $|\pi_{z,n}(\zeta)|\leqs C|c\rho(z)|^{\frac12}.
$
Now, if $\xi\in\cc$ is such that $\xi^k=\pi_{z,1}(\zeta)$, if $C$ is big enough, we have uniformly with respect to $\zeta$, $C$ and $c$
\begin{align*}
 |\pi_{z,n}(\zeta)-\pi_{z,1}(\zeta)\varphi_z(\xi,\pi_{z,2}(\zeta),\ldots,\pi_{z,n-1}(\zeta))|\leqs C|c\rho(z)|^{\frac12}\leq C^2|c\rho(z)|^{\frac12}.
\end{align*}
Moreover, since $w_z$ is a tangent vector to $bD_{\rho(z)}$ at $z$, provided $C>0$ is big enough, the point $\zeta +(\pi_{z,1}(\zeta)\varphi_z(\xi,\pi_{z,2}(\zeta),\ldots,\pi_{z,n-1}(\zeta))-\pi_{z,n}(\zeta))w_z$ belongs to $X\cap \kb{cC^4|\rho(z)|} z$.\qed
\par\medskip
Now we are interested in the case where $z$ is far from $\sing(X)$. 
\begin{proposition}\mlabel{prop2.4}
There exists $C>0$ sufficiently big  such that for all small $c>0$, for all $z$ near $0$ with $d(z,\sing(X))\geq 10 C |c\rho(z)|^{\frac12}$ and $X\cap \kbz z\neq\emptyset$, for all $\zeta\in \kbz z$ we have
 \begin{enumerate}[(i)]
  \item $|\pi_{z,1}(z)|\geq 2 C|c\rho(z)|^{\frac12}$,
  \item $|\pi_{z,1}(\zeta)|\geq  C|c\rho(z)|^{\frac12}$,
  \item $\pi_{z,1}(\zeta)$ belongs to $D(\pi_{z,1}(z),  C|c\rho(z)|^{\frac12})$.
 \end{enumerate}

\end{proposition}
\pr Firstly, for all $z$ and all $\zeta $ belonging to $\kbz z$, we have $|\pi_{z,1}(\zeta)-\pi_{z,1}(z)|\leq  C |c\rho(z)|^{\frac12}$ for some big $C>0$ which does not depend on $\zeta $ or $z$ and so $\pi_{z,1}(\zeta)$ belongs to $D(\pi_{z,1}(z), C|c\rho(z)|^{\frac12})$.

If $z$ is a point such that $d(z,\sing(X))\geq  10C |c\rho(z)|^{\frac12}$ and $X\cap \kbz z\neq\emptyset$ we prove that $|\pi_{z,1}(z)|\geq 2 C|c\rho(z)|^{\frac12}$ so that 
$|\pi_{z,1}(\zeta)|\geq  C |c\rho(z)|^{\frac12}$ for all $\zeta\in \kbz z.$

Let $\tilde \zeta $ be a point in $\kbz z\cap X$. Since $\pi_z^{-1}(0,\pi_{z,2}(\tilde \zeta),\ldots,\pi_{z,n-1}(\tilde \zeta),0)$ belongs to $\sing(X)$ and since $d(z,\sing(X))\geq 10  C|c\rho(z)|^{\frac12}$, provided $C$ is big enough, it comes 
\begin{align*}
10 C|c\rho(z)|^\frac12
&\leq |z-\pi_z^{-1}(0,\pi_{z,2}(\tilde \zeta),\ldots,\pi_{z,n-1}(\tilde\zeta),0)|\\
&\leq  \left|\pi_z^{-1}(\pi_{z,1}(\tilde \zeta),0,\ldots,0,\pi_{z,n}(\tilde \zeta))\right|+|z- \tilde\zeta|\\
&\leq \sqrt{|\pi_{z,1}(\tilde\zeta)+w_{z,1}\pi_{z,n}(\tilde\zeta)|^2+|w_{z,n}\pi_{z,n}(\tilde\zeta)|^2}+ C|c\rho(z)|^{\frac12}.
\end{align*}
Therefor
\begin{align}
\sqrt{|\pi_{z,1}(\tilde\zeta)+w_{z,1}\pi_{z,n}(\tilde\zeta)|^2+|w_{z,n}\pi_{z,n}(\tilde\zeta)|^2}\geq 9  C |c\rho(z)|^{\frac12}>0.
\label{star}
\end{align}
Since $\tilde\zeta$ belongs to $X$, if $z$ is close enough to $0$ which implies that $\tilde\zeta $ is close to $0$, Lemma \ref{lemme2.1} implies that $|\pi_{z,n}(\tilde\zeta)|\leq|\pi_{z,1}(\tilde\zeta)|$. Therefor, with (\ref{star}), we get $|\pi_{z,1}(\tilde\zeta)|\geq 3  C|c\rho (z)|^{\frac12}$.

Finally, 
since $\tilde \zeta $ belongs to $\kbz z$, $|\pi_{z,1}(\tilde \zeta)-\pi_{z,1}(z)|\leq   C |c\rho(z)|^{\frac12}$ which yields 
\begin{align*}
 |\pi_{z,1}(z)|&\geq |\pi_{z,1}(\tilde\zeta)|-|\pi_{z,1}(\tilde\zeta)-\pi_{z,1}(z)|\geq 2 C|c\rho (z)|^{\frac12}.\qed
\end{align*}

We notice that, under the hypothesis of Proposition \ref{prop2.4}, there exist holomorphic $k$-roots of $\pi_{z,1}$ on $\kbz z$. We also point out that the constant $C$ and $c$ in Propositions \ref{prop2.3} and \ref{prop2.4}  are independent. This is of importance because we will have to chose $c$ even smaller.

\begin{corollary}\label{cor2.5}
Under the hypothesis of Proposition \ref{prop2.4}, if $\xi\mapsto\xi^{\frac1k}$ is a holomorphic $k$-th root on  $D(\pi_{z,1}(z),C|c\rho(z)|^{\frac12})$, for  all $j$ we have :
 \begin{enumerate}[(i)]
  \item  if $\zeta$ belongs to $\kbz z$, the point $\pi_z^{-1}\hskip-2pt\Big(\hskip-2pt\pi_{z,1}(\zeta),\ldots,\pi_{z,n-1}(\zeta),
  \pi_{z,1}(\zeta)\varphi_z\big({\pi_{z,1}(\zeta)}^{\frac1k}\omega^j,$ $\pi_{z,2}(\zeta),\ldots,\pi_{z,n-1}(\zeta)\big)\hskip-2pt\Big)\hskip-2pt$ belongs to $X$,
  \item if  $\zeta\hskip -1pt =\hskip -1pt \pi_z^{-1}\big(\zeta'_1,\ldots,\zeta'_{n-1}, \zeta'_n)$ belongs to $\kbz z$ and if
  $|\zeta'_n-\zeta'_1\varphi_z({\zeta'_1}^{\frac1k}\omega^j,\zeta'_2,\ldots,$ $\zeta'_{n-1})|\leq 2|c\rho(z)|^{\frac12}$, then for all $\tilde\zeta$ belonging to $\kbz z$, we have, uniformly with respect to $\zeta$, $\tilde \zeta $ and $z$,
  $$|\pi_{z,1}(\tilde\zeta)\varphi_z({\pi_{z,1}(\tilde \zeta)}^{\frac1k}\omega^j,\pi_{z,2}(\tilde\zeta),\ldots,\pi_{z,n-1}(\tilde\zeta))-\pi_{z,n}(\tilde \zeta)|\leq C^2 |c\rho(z)|^{\frac12}.$$ 
  In particular, $\tilde \zeta+\left(\pi_{z,1}(\tilde\zeta)\varphi_z({\pi_{z,1}(\tilde \zeta)}^{\frac1k}\omega^j,\pi_{z,2}(\tilde \zeta),\ldots,\pi_{z,n-1}(\tilde \zeta))-\pi_{z,n}(\tilde \zeta)\right)w_z
  $ belongs to $X\cap\kb{cC^4|\rho(z)|}z$ if $C$ is big enough.
 \end{enumerate} 
\end{corollary}
\pr For $\zeta\in\kbz z$, we set $\alpha_j(\zeta )=\pi_{z,1}(\zeta)\varphi_z({\pi_{z,1}(\zeta)}^{\frac1k}\omega^j,\pi_{z,2}(\zeta),\ldots,$ $\pi_{z,n-1}(\zeta))$ as a shortcut.
By Lemma \ref{lemme2.1}, $\pi^{-1}_z(\pi_{z,1}(\zeta),\ldots,\pi_{z,n-1}(\zeta),\alpha_j(\zeta))$ belongs to $X$.

Provided $C$ is big enough, the first derivatives of $\alpha_j$ are uniformly bounded by $C$. So, if some point $\zeta$ 
belongs to $\kbz z$, then for all $\tilde \zeta \in \kbz z$, we have uniformly
\begin{align*}
 |\alpha_j(\tilde\zeta )-\pi_{z,n}(\tilde\zeta)|&\leq |\alpha_j(\tilde\zeta )-\alpha_j(\zeta )|+|\alpha_j(\zeta )-\pi_{z,n}(\zeta)|+|\pi_{z,n}(\zeta)-\pi_{z,n}(\tilde\zeta)|\\
 &\leqs C|\zeta -\tilde\zeta |+2|c\rho(z)|^{\frac12}\\
 &\leqs C|c\rho (z)|^{\frac12}\leq C^2|c\rho(z)|^{\frac12}
\end{align*}
provided $C$ is big enough.\\
Since $\tilde \zeta$ belongs to $\kbz z$ and since the vector $w_z$ is tangent to $bD_{\rho(z)}$, this implies that 
$\pi_z^{-1}\left(\hskip -2pt\pi_{z,1}(\tilde\zeta),\ldots , \pi_{z,n-1}(\tilde\zeta),\pi_{z,1}(\tilde\zeta) \varphi_z\left({\pi_{z,1}(\tilde \zeta)}^{\frac1k}\omega^j,\pi_{z,2}(\tilde\zeta),\ldots,\pi_{z,n-1}(\tilde\zeta)\right)\hskip -2pt\right)$ belongs to 
$\kb{cC^4|\rho(z)|}z$ if $C$ is big enough, independently of $c$, $\tilde \zeta$ or $z$.\qed

\medskip
We want to point out the following remarks which will be important when we will look for an upper bound of an extension. The two last propositions tell us in terms of the homogeneous geometry induced by $bD$, that if $z$ is close from $\sing(X)$, then every branches of $X$ is near $z$ and if $z$ is far from $\sing(X)$ either a branch is always far from $z$ or always close from $z$. Now we  define the current in small Koranyi balls.

\subsection{A very local definition of the current}\mlabel{verylocaldefinitionofthecurrent}

Now, for $z$ in a neighborhood $\cu (0)$ of the origin,  we construct a current $T_z$ supported in $\kbz z$ and such that $fT_z=1$.
Let $\cu(0)$ be a neighborhood of the origin such that the Propositions \ref{prop2.3} and \ref{prop2.4} and Corollary \ref{cor2.5}  hold true for all $z$ belonging to $\cu(0)$.

When $d(z,\sing(X))\geq 10C|c\rho(z)|^{\frac12}$, we choose  a holomorphic $k$-th root  $\xi\mapsto\xi^{\frac1k}$ on  $D(\pi_{z,1}(z),C|c\rho(z)|^{\frac12})$. If $\zeta'$ belongs to $\pi_z\big(\kbz z\big)$, then, by Proposition \ref{prop2.4}, $\zeta'_1$ belongs to $D(\pi_{z,1}(z),C|c\rho(z)|^{\frac12})$, so ${\zeta'_1}^{\frac1k}$ is well defined and we can set
\begin{align*}
J_z&=\left\{
j\in\{0,\ldots,k\}, \exists \zeta' \in\pi_z\big(\kbz z\big),\ \big|\zeta'_1\varphi_z({\zeta'_1}^{\frac1k}\omega^j,\zeta'_2,\ldots,\zeta'_{n-1})-\zeta'_n\big|\leq 2|c\rho(z)|^{\frac12}
\right\}. 
\end{align*}
We also denote by $\#J_z$ the cardinal of $J_z$. We use this set $J_z$ in order to pick up the branches of $X$ which are close to $z$. According to Corollary \ref{cor2.5}, for all $j\in J_z$ and all $\zeta=\pi^{-1}_z( \zeta')\in \kbz z$, we have $\big|\zeta'_1\varphi_z({\zeta'_1}^{\frac1k}\omega^j,\zeta'_2,\ldots,\zeta'_{n-1})-\zeta'_n\big|\leqs |\rho(z)|^{\frac12}$. If $j$ does not belong to $J_z$, then for all $\zeta\in \kbz z$, we have $\big|\zeta'_1\varphi_z({\zeta'_1}^{\frac1k}\omega^j,\zeta'_2,\ldots,\zeta'_{n-1})-\zeta'_n\big|\geqs |\rho(z)|^{\frac12}$

When $d(z,\sing(X))< 10C|c\rho(z)|^{\frac12}$, every branch of $X$ is close from $z$. We put  $J_z=\{0,\ldots,k-1\}$.
We then define for a smooth $(n,n)$-form $g$ supported in $\kbz z$ :
\begin{align*}
&\langle T_z,g\rangle\\
&=\frac1{\#J_z!} \int_{\zeta'\in \pi_z\left(\kbz z\right)}\hskip -2pt  \frac{\overline{\prod_{j\in J_z} \left(\zeta'_n-\zeta'_1\varphi_z({\zeta'_1}^{\frac1k}\omega^j,\zeta'_2,\ldots,\zeta'_{n-1})\right)}}{f(\pi_z^{-1}(\zeta'))}\diffp{^{\#J_z}}{\overline{\zeta'_n}^{\#J_z}}\left({\pi_z^{-1}}^*g\right)\left(\zeta'\right). 
\end{align*}
We point out that when $d(z,\sing(X))\leq 10C|c\rho(z)|^{\frac12}$, there is no holomorphic root on $\kbz z$ but $J_z=\{0,\ldots, k-1\}$ so ${\zeta'_1}^{\frac1k}$ simply needs to be any complex number such that $\left({\zeta'_1}^{\frac1k}\right)^k=\zeta'_1$.\\
Moreover, since $\prod_{j\in J_z} \left(\zeta'_n-\zeta'_1\varphi_z\left({\zeta'_1}^{\frac1k}\omega^j,\zeta'_2,\ldots,\zeta'_{n-1}\right)\right)$ is polynomial of degree $\#J_z$ in $\zeta'_n$, integrating $\#J_z$ times by parts leads to $\langle f T_z,g\rangle=\int_{\kbz z} g(\zeta )$.\\
Finally, we also notice that if $J_z=\emptyset$, the preceding definition still makes sense. In particular, if we are given a point $z$ such that $X\cap \cu(0)=\emptyset$, setting $J_z=\emptyset$ for all $z\in\cu(0)$, the previous definition gives the current 
\begin{align*}
\langle T_z,g\rangle&=\int_{\zeta\in \kbz z}  \frac1{f(\zeta)}g\left(\zeta\right). 
\end{align*}
Now we cover $\cu (0)$ with Koranyi balls and glue together the currents defined on each ball. We must define carefully the covering as in \cite{AM}.

\subsection{Koranyi covering} The Koranyi balls give to the boundary of $D$ a structure of homogeneous space. 
For $z\in D$, $v$ a unit vector in $\cc^n$, and $\varepsilon$ a positive real number, we set 
$\tau(z,v,\varepsilon)=\sup\{\tau>0,\ |\rho(z+\lambda v)-\rho(z)|<\varepsilon\text{ for all } \lambda\in\cc,\ |\lambda|<\tau\}.$
The Koranyi balls have the following homogeneous properties~:
\begin{proposition}\mlabel{propII.0.1} There exists a neighborhood $\cu$ of $bD$, a sufficiently small $c>0$ and a sufficiently big $C>0$ such that 
\begin{enumerate}[(i)]
 \item  for all $\zeta\in\cu\cap D$, ${\cal P}_{4c|\rho(\zeta)|}(\zeta)$ is included in $D$,
 \item  for all $\varepsilon>0$, all $\zeta,z\in\cu$, $\p\varepsilon\zeta\cap\p\varepsilon z\neq
\emptyset$ implies $\p\varepsilon z\subset \p{C\varepsilon}\zeta$,
\item  for all $\varepsilon>0$ sufficiently small, all
$z\in\cu$, all $\zeta\in \p{\varepsilon}z$ we have uniformly $|\rho(z)-\rho(\zeta )|\leqs
 \varepsilon$,
\item for all $\varepsilon >0$, all unit vector $v\in\cc^n$, all $z\in\cu$ and all $\zeta\in\p \varepsilon z$, $\tau(z,v,\varepsilon )
\eqs\tau(\zeta,v,\varepsilon )$ uniformly with respect to $\varepsilon,$ $z$ and $\zeta $.
\end{enumerate}
\end{proposition}
For $\cu$ given by Proposition \ref{propII.0.1} and $z$ and $\zeta$ belonging to $\cal U$, we set $\delta(z,\zeta)=\inf\{\varepsilon>0,
\zeta\in \p\varepsilon z\}$. Proposition \ref{propII.0.1} implies  that $\delta$ is a pseudo-distance in the
following sense: 
\begin{proposition}\label{propII.0.2}
For $\cu$ and $C$ given by Proposition \ref{propII.0.1} and for all $z,\ \zeta$ and $\xi$ belonging to $\cu$  we have
$$\frac1{C}\delta(\zeta,z)\leq \delta(z,\zeta)\leq C \delta(\zeta,z)$$
 and
$$\delta(z,\zeta)\leq C(\delta(z,\xi)+\delta(\xi,\zeta)).$$
\end{proposition}

Let $\cu$ be a subset of $\cc^n$ and $\varepsilon_0$ be a small positive number. We cover $\cu\cap (D\setminus D_{-\varepsilon_0})$ with a family of \ko balls 
$\kbc$ where $c$ is a small positive real number. We  assume that $c$ is so small that for all $z$ and all $\zeta\in\kbz z$ we have $\frac12|\rho(z)|<|\rho(\zeta)|<\frac32|\rho(z)|$. This construction uses classical ideas
of the theory of homogeneous spaces and is analogous to the construction of the covering of
\cite{BCD}.\\
Let $\varepsilon_0$, $\kappa$ and $c$ be positive real numbers sufficiently small. We construct a sequence of points of $\cu\cap(D\setminus D_{-\varepsilon _0})$
as follows.\\
Let $k$ be a non-negative integer and choose arbitrarily $z_1^{(k)}$ in $bD_{-(1-c\kappa )^k\varepsilon _0}$.\\
When $z_1^{(k)},\ldots, z_j^{(k)}$ are chosen, there are two possibilities. Either for all $z\in bD_{-(1-c\kappa )^k\varepsilon _0}\cap \cu$ there exists
$i\leq j$ such that $\delta(z,z_i^{(k)})<c\kappa (1-c\kappa )^k\varepsilon _0$ and the process ends here, or there exists
$z\in bD_{-(1-c\kappa )^k\varepsilon _0}\cap \cu$ such that  for all $i\leq j$ we have $\delta(z,z_i^{(k)})\geq c\kappa (1-c\kappa )^k\varepsilon _0$ 
and we chose $z^{(k)}_{j+1}$ among these points. Since $bD_{-(1-c\kappa )^k\varepsilon _0}\cap\cu$ is bounded, this process stops at some rank $n_k$.\\
We thus have constructed a sequence $(z_j^{(k)})_{k\in\nn, j\in\{1,\ldots,n_k\}}$ such that
\begin{enumerate}[(i)]
 \item \label{seqi} For all $k\in\nn$, and all $j\in\{1,\ldots, n_k\}$, $z_j^{(k)}$ belongs to $bD_{-(1-c\kappa )^k\varepsilon _0}\cap \cu$.
 \item \label{seqii} For all $k\in\nn$, all $i,j\in \{1,\ldots, n_k\}$, $i\neq j$, we have $\delta(z_i^{(k)},z_j^{(k)})\geq c\kappa 
 (1-c\kappa )^k\varepsilon _0$.
 \item \label{seqiii} For all $k\in\nn$, all $z\in bD_{-(1-c\kappa )^k\varepsilon _0}$, there exists $j\in\{1,\ldots, n_k\}$ such that 
 $\delta(z,z_j^{(k)})<c\kappa (1-c\kappa )^k\varepsilon _0$.
\end{enumerate}
For such sequences, we prove the following proposition.
\begin{proposition}\mlabel{propmax}
For $\kappa >0$ and $c>0$ small enough, let $\left(z_j^{(k)}\right)_{k\in\nn,j\in\{1,\ldots, n_k\}}$ be a sequence which satisfies (\ref{seqi}),
(\ref{seqii}) and (\ref{seqiii}). Then
\begin{enumerate}[(a)]
 \item \label{propmax1} $D\setminus D_{-\varepsilon _0}\cap\cu$ is included in $ \cup_{k=0}^{+\infty} \cup_{j=1}^{n_k} 
 \kb{c|\rho (z_j^{(k)})|} {z_{j}^{(k)}}$,
 \item \label{propmax2} there exists $M\in\nn$ such that for all $z\in D\setminus D_{-\varepsilon _0}\cap \cu$, $\kb{5cC^4 |\rho (z)|}z$ intersects at most $M$ Koranyi balls $\kb{ 5cC^4 |\rho (z_j^{(k)})|}{z_{j}^{(k)}}$.
\end{enumerate}
\end{proposition}
\pr 
We first prove that (\ref{propmax1}) holds. For $z\in D\setminus D_{\varepsilon _0}$, let $k\in\nn$ be such that 
$$(1-c\kappa )^{k+1}\varepsilon _0\leq|\rho (z)|<(1-c\kappa )^k\varepsilon _0,$$
and let $\lambda\in\cc$ be such that $\zeta=z+\lambda \eta_z$ belongs to $bD_{-(1-c\kappa)^k\varepsilon_0}$. 
On the one hand, the assumption $(\ref{seqiii})$ implies that there exists $j\in\{1,\ldots, n_k\}$ such that $\delta\left(\zeta,z_j^{(k)}\right)
\leq c\kappa (1-c\kappa )^k\varepsilon _0$. On the other hand we have $|\lambda|=\delta(z,\zeta)\leq \tilde C c\kappa(1-c\kappa)^k\varepsilon_0$ where
$\tilde C$  neither depends on $z$ nor on $\zeta $ nor on $c\kappa$. These two inequalities yield
\begin{eqnarray*}
 \delta\left(z,z^{(k)}_j\right)&\leq&C\left(\delta(z,\zeta )+\delta(\zeta,z^{(k)}_j)\right)\\
&\leq& \kappa cC(1-c\kappa )^k\varepsilon _0 (\tilde C +1)\\
&\leq& c \left|\rho \left(z_j^{(k)}\right)\right|
\end{eqnarray*}
provided $\kappa$ is small enough. Therefore $z$ belongs to $\kb{c|\rho (z_j^{(k)})|}{z_j^{(k)}}$ and (\ref{propmax1}) holds.\\
We now prove (\ref{propmax2}). Let $z$ be a point of $D\setminus D_{-\varepsilon _0}\cap\cu(0)$. For all $\zeta \in{\cal P}_{5cC^4|\rho(z)|}(z)$, if $c$ is small enough, proposition \ref{propII.0.1} yields
$$\frac12 |\rho (z)|\leq |\rho (\zeta )|\leq 2|\rho (z)|.$$
The same inequalities hold for all $z^{(k)}_j$ and all $\zeta\in \kb{5cC^4 |\rho (z_j^{(k)})|}{z_j^{(k)}}$. 
Thus if the intersection $\kb{5cC^4|\rho (z_j^{(k)})|}{z_j^{(k)}}\cap \kb{5cC^4|\rho (z)|}z$ is not empty, we have
$$\frac14|\rho (z)| \leq (1-c\kappa )^k\varepsilon_0\leq 4 |\rho (z)|.$$
Therefore $\frac{\ln\left(\frac{\rho(z)}{\varepsilon_0}\right)-\ln 4}{\ln(1-c\kappa)}\geq k\geq \frac{\ln\left(\frac{\rho(z)}{\varepsilon_0}\right)+\ln 4}{\ln(1-c\kappa)}$ so $k$ can take only a finite number of values.\\
For such a $k$, we set $I_k=\left\{j\in\{1,\ldots,n_k\},\ \kb{5cC^4 |\rho (z_j^{(k)})|}{z_j^{(k)}}\cap \kb{5cC^4
|\rho (z)|}z\neq \emptyset\right\}$. Assertion (\ref{propmax2}) will be proved provided we show that $\#I_k$, the cardinal of $I_k$, is bounded
uniformly with respect to $k$ and $z$.\\
We denote by $\sigma$ the area measure on $bD_{-(1-c\kappa )^k\varepsilon _0}$. For all $i,j\in I_k$ distinct, we have
$ \delta\left(z_{i}^{(k)},z_j^{(k)}\right)\geq c\kappa (1-c\kappa )^k\varepsilon _0$. So, provided $C$ is big enough,  we have
\begin{align*}
\sigma \left(\cup_{j\in I_k}\kb{5cC^4\left|\rho \left(z_j^{(k)}\right)\right|}{z_j^{(k)}}\cap 
bD_{-(1-c\kappa )^k\varepsilon _0}\right)
& \geq\sigma \left(\cup_{j\in I_k}\kb{\frac {c\kappa}{C} (1-c\kappa )^k\varepsilon _0}{z_j^{(k)}}\cap bD_{-(1-c\kappa )^k
\varepsilon _0}\right)\\
& \geq \# I_k \cdot\left(\frac {c\kappa}{C}  (1-c\kappa )^k\varepsilon _0\right)^n.
\end{align*}
Now we look for an upper bound of $\sigma \left(\cup_{j\in I_k}\kb{5cC^4|\rho (z_j^{(k)})|}{z_j^{(k)}}\cap bD_{-(1-c\kappa )^k\varepsilon_0}
\right)$.
We  fix $j_0\in I_k$. For all $j\in I_k,$ since $\kb{5cC^4|\rho (z_j^{(k)})|}{z_j^{(k)}}\cap \kb{5cC^4 |\rho (z)|}z$
and $\kb{5cC^4 |\rho (z_{j_0}^{(k)})|}{z_{j_0}^{(k)}}\cap \kb{5cC^4|\rho (z)|}z$ are not empty, we have 
\begin{align*}
\delta\left(z_{j_0}^{(k)},z_j^{(k)}\right)
&\leqs \delta\left(z_{j_0}^{(k)},z\right) +\delta\left(z,z_j^{(k)}\right)\\
&\leqs 5cC^4 \left(\left|\rho \left(z_{j_0}^{(k)}\right)\right|+\left|\rho \left(z_{j}^{(k)}\right)\right|\right)\\
&\leqs  5cC^4(1-c\kappa )^k\varepsilon _0
\end{align*}
uniformly with respect to $k$, $j$ and $j_0$. Thus there exists $K$ neither depending on $z$, nor on $j$, nor on  $j_0$ nor on $k$ such that 
$\kb{5cC^4|\rho (z_j^{(k)})|}{z_j^{(k)}}
\subset \kb{5cC^4K |\rho (z_{j_0}^{(k)})|}{z_{j_0}^{(k)}}$. Therefore 
\begin{align*}
\sigma \left(\cup_{j\in I_k}\kb{5cC^4|\rho (z_j^{(k)})|}{z_j^{(k)}}\cap bD_{-(1-c\kappa )^k\varepsilon _0}\right) 
&\leq
\sigma \left({\cal P}_{5cC^4K|\rho (z_{j_0}^{(k)})|}{z_{j_0}^{(k)}}\cap bD_{-(1-c\kappa )^k\varepsilon _0}\right)\\
&\leqs \left(KcC^4(1-c\kappa )^k\varepsilon _0\right)^n
\end{align*}
which yields $\#I_k\leqs \left(\frac{C^5}\kappa\right)^{n}$.\qed
\par\medskip
The covering property (\ref{propmax1}) allows us to settle the following definition :
\begin{definition}
 Let $\cu$ be any subset of $\cc^n$. If the sequence $(z_j)_{j\in\nn}$ can be renumbered such that (\ref{seqi}),  (\ref{seqii}) and (\ref{seqiii}) hold true, the family $\left(\kb{c|\rho (z_j)|}{z_j}
 \right)_{j\in\nn}$ will be called a $c$-covering of ${\cu}\cap (D\setminus D_{-\varepsilon_0})$.
\end{definition}
If we are given an open $\cu$ and a $c$-covering $\left(\kb{c|\rho (z_j)|}{z_j}
 \right)_{j\in\nn}$ of $\cu\cap(D\setminus D_{\varepsilon_0})$, we will need to know how many Koranyi balls of a given diameter cover it.  In this spirit, we prove the following lemma.
\begin{lemma}\label{nombre_boules}
 Let $\cu$ be an open set, let $\left(\kb{c|\rho (z_j)|}{z_j}  \right)_{j\in\nn}$
 be a $c$-covering of $\cu\cap(D\setminus D_{\varepsilon_0})$ and let $z$ be a point in $\cu\cap(D\setminus D_{\varepsilon_0})$.
 Let us renumber the point $(z_j)_{j\in\nn}$ in the following way :\\
 Let $j_0$ be an integer such that $(1-c\kappa)^{-j_0}\varepsilon_0\leq|\rho(z)|\leq (1-c\kappa)^{-j_0-1}\varepsilon_0$, let $z_1^{i,j},\ldots,z^{i,j}_{m_{i,j}}$, $i\in\nn$, $j\in\zz$, be the  points of the covering such that
 \begin{itemize}
  \item $\rho(z_m^{i,j})=-(1-c\kappa)^{j-j_0}\varepsilon_0$,
  \item $\delta(z_m^{i,j},z)$ belongs to $[ic(1-c\kappa)^{j-j_0}\varepsilon_0,(1+i)c(1-c\kappa)^{j-j_0}\varepsilon_0[,$
  \item $\delta(z_m^{i,j},z)\leq \varepsilon_0$.
 \end{itemize}
For $j\geq j_0$, let 
$i_0(j)$ be the non-negative integer such that $i_0(j)\kappa(1-c\kappa)^{j-j_0}<1
\leq (1+i_0(j))\kappa(1-c\kappa)^{j-j_0}$.\\
Then 
\begin{enumerate}[(i)]
 \item \mlabel{premierpoint} $\kb{\frac{\varepsilon_0}{2C}}z\cap D\subset\cup_{j=j_0}^{+\infty}\cup_{i=0}^{i_0(j)}\cup_{m=1}^{m_{i,j}}\kb{c|\rho(z^{i,j}_m)|} {z_m^{i,j}}$,
 \item \mlabel{secondpoint} $m_{i,j}\leqs i^n$ uniformly with respect to $z_0,z,i$ and $j$.
\end{enumerate}
\end{lemma}
\pr Point (\ref{premierpoint}) can be shown exactly as lemma 4.2 in \cite{AM} because this property relies only on the homogeneous properties of the Koranyi balls and not on the dimension. For the second point (\ref{secondpoint}), the only difference with \cite{AM} is that $\sigma(bD_\varepsilon\cap \kb r z)$, the area of $bD_\varepsilon\cap \kb r z$, if non empty, is of order $r\cdot\left({\sqrt{r}}\right)^{2(n-1)}=r^n$ and so we get a power $n$ instead of $2$ in dimension $2$.\qed 

\subsection{The global definition of the current}
For any point $p$ of $\overline{D}$, we first define locally currents $T^{(\cu)}$ supported on a neighborhood $\cu$ of $p$ and such that $fT^{(\cu)}=1$.
\par\bigskip
\begin{trivlist}
\item[\bf First case :] If we are given a point $p$ of $bD\cap X$, let $\cu$ be a neighborhood of $p$ such that the Propositions \ref{prop2.3}, \ref{prop2.4} and \ref{propII.0.1} hold true for all $z$ in $\cu$. We fix an arbitrary $c$-covering $\left(\kb{c|\rho (z_j)|}{z_j} \right)_{j\in\nn}$ of ${\cu}\cap D$ and a partition of unity $(\psi_j)_{j\in\nn}$ associated to this $c$-covering. We mention here that any partition of unity associated with the covering $\left(\kb{c|\rho (z_j)|}{z_j} \right)_{j\in\nn}$ will give a current $T^{(\cu)}$ such that $fT^{(\cu)}=1$ on ${\cu}\cap (D\setminus D_{-\varepsilon_0})$. However, in order to have a good extension, we must choose a partition of unity such that for all $j$, if $w_1=\eta_{z_j}, w_2,\ldots, w_n$ is an orthonormal basis, then for all multi-index $\alpha$ and $\beta$, $\left|\diffp{^{|\alpha|+|\beta|} \psi_j}{w^\alpha\partial\overline{w}^\beta}(z)\right|\leqs \frac1{|\rho(z_j)|^{\alpha_1+\beta_1+\frac{\alpha_2+\beta_2+\ldots+\alpha_n+\beta_n}2}}$. We then set
$$T^{(\cu)}=\sum_{j\in\nn} \psi_j T_{z_j},$$
where $T_{z_j}$ is the current defined in Subsection \ref{verylocaldefinitionofthecurrent}. Since $fT_{z_j}=1$ for all $j$, we trivially have $fT^{(\cu)}=1$ on $\cu\cap D$. 
Moreover, if the support of a form $g$ does not intersect $\overline D$, it intersects none of the Koranyi ball $\kb {c|\rho(z_j)|}{z_j}$ and we have $\langle T^{(\cu)},g\rangle=0$. Therefor the support of $T^{(\cu)}$ is included in $\overline D$.
\item[\bf Second case :] 
If we are given a point $p$ of $bD\setminus X$, let $\cu$ be a neighborhood of $p$ such that $\cu\cap X=\emptyset$. We then define for a smooth $(n,n)$-form $g$ supported in $\cu$ :
\begin{align*}
\langle T^{(\cu)},g\rangle&=\int_{\zeta\in\cu\cap D} \frac1{f(\zeta)} g(\zeta). 
\end{align*}
In this case we trivially have $fT^{(\cu)}=1$ on $\cu\cap D$ and $T^{(\cu)}$ is trivially supported in $\overline D$.
\par\medskip
By compacity of $bD$, we get finitely many open $\cu_1,\ldots,\cu_{N}$ given by the first and second cases such that $bD$ is included in $\bigcup_{j=1}^N \cu_j$. There exists a small real number that we still denote by $\varepsilon_0$ such that $D\setminus D_{-\varepsilon_0}\subset  \bigcup_{j=1}^N\cu_j$.
\item[\bf Third case :] 
For any point $p\in \overline{D_{-\varepsilon_0}},$ maybe after a linear change of coordinates, there is a neighborhood $\cu$ of $p$ included in $D$, a non vanishing holomorphic function $u$
and
a  Weierstrass polynomial $P_{\cu}(z)=z^l_n+a_{1}(z_1,\ldots,z_{n-1})z_n^{l-1}+\ldots+ a_{l}(z_1,\ldots,z_{n-1})$
such that $f=uP_\cu$ on $\cu$.
Now for a smooth $(n,n)$-form $g$ supported in $\cu$, we put :
\begin{align*}
\langle T^{(\cu)},g\rangle&=\frac1{l!}\int_{\zeta\in\cu} \frac{\overline{P_\cu(\zeta)}}{f(\zeta)}\diffp{^l g}{\overline{\zeta_n}^l}(\zeta). 
\end{align*}
As in the very local definition of the current in Subsection \ref{verylocaldefinitionofthecurrent}, after $l$ integrations by parts, we get $fT^{(\cu)}=1$. Moreover, $T^{(\cu)}$ is again supported in $\overline D$.

Finally, we choose a finite covering of $\overline{D_{-\varepsilon_0}}$ and putting it together with the finite covering of $bD$, we get a covering of $\overline D$ that we still denote by $\cu_1,\ldots,\cu_N$ and the associated currents $T^{(\cu_1)},\ldots,T^{(\cu_N)}$. Using a partition of unity $(\chi_j)_{j=1,\ldots, N}$ associated with this covering, we set 
$$T=\sum_{j=1}^N\chi_j T^{(\cu_j)}.$$
\end{trivlist}

We have here defined a current such that $fT=1$ on $D$. We will need to apply $T$ to forms which are not necessarily $C^\infty$ smooth and with support included in $\overline D$ and not only in $D$. We now prove that $T$ can be applied to forms whose support is included in $\overline D$ and of class $C^{k+1}$ where $k$ is the maximal order of the singularities of $X$.
\par\medskip
Let $\cu$ be an open set among the $\cu_j$ which define $T$ and let $g$ be a $(n,n)$-form of class $C^{k+1}$ supported in $\overline D\cap \cu$. 
\par\medskip
If we are given an open set $\cu$ as in the second case, then $\frac1{f}$ is bounded and so
\begin{align*}
 \left|\int_{\cu} \frac1{f(\zeta)} g(\zeta)\right|&\leqs \|g\|_{C^{k+1}},\end{align*}
where $\left\|g \right\|_{C^{k+1}}$ is the usual $C^{k+1}$ norm of $g$.
\par\medskip
If $\cu$ is as in the third case, then $\frac{\overline{P_\cu}}{f}$ is bounded and 
\begin{align*}
 \left|\int_{\zeta\in\cu} \frac{\overline{P_\cu(\zeta)}}{f(\zeta)}\diffp{^l g}{\overline{\zeta_n}^l}(\zeta)\right|
 &\leqs \left\|g \right\|_{C^{k+1}}.
 \end{align*}
\par\medskip
Now if $\cu$ is as in the first case, without restriction, we may assume that $\cu$ is of the form $\kb{\frac{\varepsilon_0}{2C}}z.$ Let $\left(\kb{c|\rho (z_j)|}{z_j}  \right)_{j\in\nn}$ be the $c$-covering  and let $(\psi_j)_{j\in\nn}$ be the partition of unity used in the definition of $T^{(\cu)}$. 
Then for all $j$, the quotient 
$\frac{\overline{\prod_{l\in J_{z_j}} \left(\zeta'_n-\zeta'_1\varphi_{z_j}({\zeta'_1}^{\frac1k}\omega^l,\zeta'_2,\ldots,\zeta'_{n-1})\right)}}{f(\pi_{z_j}^{-1}(\zeta'))}$ is bounded by $|\rho(z_j)|^{-\frac {k-\# J_{z_j}}2}$.\\
On the other hand, since $g$ vanishes outside of $D$ and since $bD$ is smooth (at least of classe $C^{k+1}$), for all $z\in D$ and all linear differential operator $\Delta_l$ of order $l\leq k$, $|\Delta_l g(z)|\leqs |\rho(z)|^{k+1-l}$, uniformly with respect to $z$. Proposition \ref{propII.0.1} implies if $c$ is sufficiently small that $|\rho(z)|\eqs|\rho(z_j)|$ uniformly with respect to $z$ in $\kb{c|\rho (z_j)|}{z_j}$. 
Moreover, $\left|\diffp{^l\psi_j}{w_{z_j}^l}(\zeta)\right|\leqs |\rho(z_j)|^{-\frac l2}$ for all $l$.
Therefore
\begin{align*}
&\left| \int_{\zeta'\in \pi_{z_j}\left(\kbz {z_j}\right)}\hskip -2pt  \frac{\overline{\prod_{l\in J_{z_j}} \left(\zeta'_n-\zeta'_1\varphi_{z_j}({\zeta'_1}^{\frac1k}\omega^j,\zeta'_2,\ldots,\zeta'_{n-1})\right)}}{f(\pi_{z_j}^{-1}(\zeta'))}\diffp{^{\#J_{z_j}}}{\overline{\zeta'_n}^{\#J_{z_j}}}\left({\pi_{z_j}^{-1}}^*\varphi_{z_j}g\right)\left(\zeta'\right)\right|\\
&\leqs \int_{\kbz {z_j}}  |\rho(z_j)|^{-\frac{k-\#J_{z_j}}2} |\rho(z_j)|^{k+1-\#J_{z_j}} \|g\|_{C^{k+1}}dV(\zeta)\\
&\leqs |\rho(z_j)|^{n+2} \|g\|_{C^{k+1}}.
\end{align*}
Now, if we renumber the sequence $(z_j)_{j\in\nn}$ as in Lemma \ref{nombre_boules}, we get
\begin{align*}
 &\sum_{j=0}^{+\infty}\left| \int_{\zeta'\in \pi_{z_j}\left(\kbz {z_j}\right)}\hskip -4pt  \frac{\overline{\prod_{l\in J_{z_j}} \left(\zeta'_n-\zeta'_1\varphi_{z_j}({\zeta'_1}^{\frac1k}\omega^j,\zeta'_2,\ldots,\zeta'_{n-1})\right)}}{f(\pi_{z_j}^{-1}(\zeta'))}\diffp{^{\#J_{z_j}}}{\overline{\zeta'_n}^{\#J_{z_j}}}\left({\pi_{z_j}^{-1}}^*\varphi_{z_j}g\right)\left(\zeta'\right)\right|\\
 &\leqs \sum_{j=j_0}^{+\infty}\sum_{i=0}^{i_0(j)}\sum_{m=1}^{m_{i,j}} |\rho(z_m^{i,j})|^{n+2} \|g\|_{C^{k+1}}\\
 &\leqs \|g\|_{C^{k+1}}\sum_{j=j_0}^{+\infty}\sum_{i=0}^{i_0(j)}\sum_{m=1}^{m_{i,j}} \left((1-c\kappa)^{j-j_0}\varepsilon_0\right)^{n+2}\\
 &\leqs  \|g\|_{C^{k+1}}\varepsilon_0^{n+2}\sum_{j=j_0}^{+\infty} (1-c\kappa)^{(j-j_0)(n+2)} \sum_{i=0}^{i_0(j)}i^{n}\\
 &\leqs  \|g\|_{C^{k+1}}\varepsilon_0^{n+2}\sum_{j=j_0}^{+\infty} (1-c\kappa)^{(j-j_0)(n+2)} i_0(j)^{n+1}\\
 &\leqs \|g\|_{C^{k+1}}\varepsilon_0^{n+2}\sum_{j=j_0}^{+\infty} (1-c\kappa)^{j-j_0}\\
 &\leqs \|g\|_{C^{k+1}}.
\end{align*}
Therefor, we can apply $T$ to forms of class $C^{k+1}$, where $k$ is the maximum order of the singularities of $X$, and which vanish outside of $D$.
\subsection{Berndtsson-Andersson reproducing kernel in $\cc^n$}
\mlabel{subsection2.5}
We now recall the definition of the Berndtsson-Andersson kernel of $D$ when $D$ is a  strictly pseudoconvex domain of $\cc^n$.

We denote by $g^{(1)}$ the support function for strictly pseudoconvex domain,  given in Theorem 16 of \cite{For} and the remarks following it, times $-\frac12$ and by $(g^{(1)}_j)_{j=1,\ldots, n}$ its Hefer decomposition. If $bD$ is of class $C^r$ , the functions $g^{(1)}(\zeta,z)$ and $g^{(1)}_j(\zeta,z)$ are of class $C^{r-1}$ in a neighborhood of $D\times D$, holomorphic with respect to $z$ when $\zeta$ is fixed and $g^{(1)}$ satisfies
\begin{enumerate}[(i)]
 \item $g^{(1)}(\zeta,\zeta)=0$ for all $\zeta$ in a neighborhood of $D$,
 \item for all $\zeta$ and $z$ in a compact set $K$ of $D$ :
$$2 \re\ g^{(1)}(\zeta,z)\leq \rho(z)- \rho(\zeta)-2\delta_K d(\zeta,z)^2,$$
where $d(\cdot,\cdot)$ is a distance on a neighborhood of $D$ and $\delta_K>0$ depends on $K$.
\end{enumerate}
 We cannot use these functions in order to define our extension operator, because  we would fail to get good estimates of our extension. We need a more explicit support function near the boundary as the Levi polynomial $F$. We set for $\zeta$ near $bD$ and $z$ in $\cc^n$ :
$$F(\zeta,z)=
-\sum_{j=1}^n\diffp{\rho}{\zeta_j}(\zeta)(z_j-\zeta_j)
+\frac12\sum_{j,k=1}^n\diffp{^2\rho}{\zeta_j\partial\zeta_k}(\zeta)(\zeta_j-z_j)(\zeta_k-z_k).$$
For all $\zeta$ and $z$ in a neighborhood of $bD$, such that $|\zeta-z|$ is sufficiently small, we have
$$2\re\ F(\zeta,z)\leq \rho(z)-\rho(\zeta) - 2\beta |\zeta-z|^2$$
where $\beta$ is strictly positive.
Since this estimate holds true only for $\zeta$ and $z$ close from each other, we have to use a global version of the Levi polynomial. From the proof of Proposition VII.3.1 of \cite{Ran}, there exist a neighborhood $\cu(bD)$ of $bD$, a neighborhood $\cu(\overline D)$ of $\overline D$, $\varepsilon>0$, a function $g^{(0)}$ defined on $\cu(bD)\times \cu(\overline D)$ and a function $v$ defined on   $\cu(\overline D)\times \cu(\overline D)$ such that
\begin{enumerate}[(i)]\setcounter{enumi}{2}
 \item $g^{(0)}(\zeta,z)$ is holomorphic with respect to $z$ in $\cu(\overline D)$,
 \item $\re\ g^{(0)}(\zeta,z)> 0$ for all $(\zeta,z)\in\cu(bD)\times\cu(\overline D)$ such that $|\zeta-z|\geq \varepsilon$,
 \item $g^{(0)}(\zeta,z)=\frac12 F(\zeta,z)\frac1{1+F(\zeta,z)(v-\sup_{\cu(\overline D)}|v|)}$ and $\re\ g^{(0)}(\zeta,z)\leq \re\ F(\zeta,z)$ for all  $(\zeta,z)\in\cu(bD)\times\cu(\overline D)$ such that $|\zeta-z|\leq \varepsilon$.
\end{enumerate}
Let $g^{(0)}_1,\ldots,g^{(0)}_n$ be the Hefer decomposition of $g^{(0)}$.
Finally, let $\chi$ be a smooth cutoff function such that $\chi=1$ near $bD$ and $\chi=0$ outside of the neighborhood $\cu(bD)$ where $g^{(0)}$ is defined. Then, we put
\begin{align*}
 g(\zeta,z)&=\chi(\zeta) g^{(0)}(\zeta,z)+(1-\chi(\zeta))g^{(1)}(\zeta,z),\\
 g_j(\zeta,z)&=\chi(\zeta) g_j^{(0)}(\zeta,z)+(1-\chi(\zeta))g_j^{(1)}(\zeta,z),\\
 G(\zeta,z)&=\frac1{\rho(\zeta)}\sum_{j=1}^n g_j(\zeta,z)d\zeta_j.
\end{align*}
Then 
\begin{itemize}
 \item $g$ is defined in a neighborhood of $\overline D\times\overline D,$
 \item There exists $\gamma>0$ such that $\re\ G(\zeta,z)+1\geq \frac1{\rho(\zeta)}\left(\frac{\rho(\zeta)+\rho(z)}2-\gamma|\zeta-z|^2\right)>0$ for all $(\zeta,z)\in \cu(D)\times \cu(D)$,
 \item $\sum_{j=1}^n g_j(\zeta,z)(\zeta_j-z_j)=g(\zeta,z)$.
\end{itemize}
We then can define
the Berndtsson-Andersson reproducing kernel by setting for
 arbitrary positive integers $N$ and $k$ and all $\zeta,z\in D$~:
$$P^{N,k}(\zeta,z)=C_{N,k} \left(\frac{\rho(\zeta)}{g(\zeta,z)+\rho(\zeta)}\right)^{N+k}\left(\overline\partial G(\zeta,z)\right)^k,$$
where $C_{N,k}\in \cc$ is a constant. We also set $P^{N,k}(\zeta,z)=0$ for all $z\in D$ and all $\zeta\notin D$. 
The following representation formula holds (see \cite{BA}):
\begin{theorem}
 For all $g\in \co(D)\cap C^\infty(\overline D)$ we have
$$g(z)=\int_D g(\zeta)P^{N,n}(\zeta,z).$$
\end{theorem}

\subsection{The extension operator}\label{subsection2.6}
We now come to the definition of the extension operator. Let $b(\zeta,z)=\sum_{j=1}^nb_j(\zeta,z)d\zeta_j$ be a holomorphic $(1,0)$-form such that $\sum_{j=1}^n b_j(\zeta,z)(z_j-\zeta_j)=f(z)-f(\zeta)$. 
As in \cite{AM}, for a function $h$ holomorphic on $X\cap D$ and admitting a smooth extension $\tilde h$ which satisfies the assumptions of Theorem \ref{th0},  we define the extension $E_N(h)$ of $h$ by setting 
$${E_N}(h)(z)=C_1 \overline\partial T[\tilde hb(\cdot,z)\wedge P^{N,n-1}(\cdot,z)],\qquad\forall z\in D,$$
where $C_1$ is a suitable constant (see \cite{Maz1}). It was shown in \cite{AM} that $E_N(h)$ is indeed an extension of $h$ and does not depend on the choice of the smooth extension $\tilde h$, provided it satisfies the hypothesis of Theorem \ref{th0}.
However, we have to make some restriction on the choice of $b$ in order to have a good extension. We need that $b$ satisfies the following lemma.
\begin{lemma}\label{lemme Hefer} For all positive integer $k$ and all function $f$, holomorphic in a pseudoconvex neighborhood of $\overline D$, there exist functions $b_1,\ldots, b_n$ holomorphic  in a neighborhood of $\overline{D}\times \overline D$ such that  
 \begin{enumerate}[(i)]
  \item $\sum_{j=1}^n b_j(\zeta,z)(z_j-\zeta_j)=f(z)-f(\zeta)$,
  \item $b_j(\zeta,z)=\sum_{|\alpha|\leq k} \frac1{\alpha!(|\alpha|+1)} \diffp{^{|\alpha|+1} f}{\zeta^\alpha\partial\zeta_j}(\zeta)(z-\zeta)^\alpha +O(|\zeta-z|^{k+1}).$
 \end{enumerate}
\end{lemma}
\pr Since $f$ is holomorphic in a pseudoconvex neighborhood of $\overline D$, there exist $n$ functions  $\tilde b_1,\ldots,\tilde b_n$ holomorphic in a neighborhood of $\overline D\times\overline D$,  such that $f(z)-f(\zeta)=\sum_{j=1}^n \tilde b_j(\zeta,z)(z_j-\zeta_j)$. Therefor, there exists $R>0$ such that for all $\zeta\in\overline D$ and all $u\in \cc^n$ with $|u|<R$ the following equalities hold :
\begin{align*}
 f(\zeta+u)-f(\zeta)&=\sum_{l=1}^n\tilde b_l(\zeta,\zeta+u)u_k\\
 &=\sum_{l=1}^n \sum_{j=1}^{+\infty}\sum_{|\alpha|=j}\frac1{\alpha!}\diffp{^j\tilde b_l}{z^\alpha}(\zeta,\zeta)u^\alpha u_l.
\end{align*}
On the other hand, for all $u\in \cc^n$ with $|u|<R$
\begin{align*}
 f(\zeta+u)-f(\zeta)&=\sum_{j=1}^{+\infty}\sum_{|\alpha|=j}\frac1{\alpha!}\diffp{^jf}{\zeta^\alpha}(\zeta)u^\alpha.
\end{align*}
These yields for any multi-index $\alpha$ with $|\alpha|>0$ :
\begin{align}
 \frac1{\alpha!}\diffp{^{|\alpha|}f}{\zeta^\alpha}(\zeta)
 &=\sum_{\over{j=1}{\alpha_j\geq 1}}^n \frac1{\alpha_1!\ldots(\alpha_j-1)!\ldots\alpha_n!} \diffp{^{|\alpha|-1}\tilde b_j}{z^{(\alpha_1,\ldots, \alpha_j-1,\ldots,\alpha_n)}}(\zeta,\zeta)\nonumber\\
 &=\sum_{{j=1}}^n \frac{\alpha_j}{\alpha!} \diffp{^{|\alpha|-1}\tilde b_j}{z^{(\alpha_1,\ldots, \alpha_j-1,\ldots,\alpha_n)}}(\zeta,\zeta).\label{eq3}
\end{align}
We now define
$b_j(\zeta,z)=\tilde b_j(\zeta,z)-\sum_{|\alpha|\leq k}\frac1{\alpha!}\diffp{^{|\alpha|}\tilde b_j}{z^\alpha}(\zeta,\zeta)(z-\zeta)^\alpha +\sum_{|\alpha|\leq k} \frac1{\alpha!} \diffp{^{|\alpha|+1}f}{\zeta^\alpha\partial\zeta_j}(\zeta) \frac{(z-\zeta)^\alpha}{|\alpha|+1}.$
We have 
\begin{align*}
b_j(\zeta,z)- \sum_{|\alpha|\leq k} \frac1{\alpha!} \diffp{^{|\alpha|+1}f}{\zeta^\alpha\partial\zeta_j}(\zeta) \frac{(z-\zeta)^\alpha}{|\alpha|+1}&= \tilde b_j(\zeta,z)-
\sum_{|\alpha|\leq k}\frac1{\alpha!}\diffp{^{|\alpha|}\tilde b_j}{z^\alpha}(\zeta,\zeta)(z-\zeta)^\alpha \\
&=O(|\zeta-z|^{k+1}).
\end{align*}
Now we check that $f(z)-f(\zeta)=\sum_{j=1}^n b_j(\zeta,z)(z_j-\zeta_j)$. The definition of $b$ implies that
\begin{align*}
 \sum_{j=1}^n b_j(\zeta,z)(z_j-\zeta_j)
 =&f(z)-f(\zeta)-\sum_{j=1}^n\sum_{|\alpha|\leq k}\frac1{\alpha!}\diffp{^{|\alpha|} \tilde b_j}{z^\alpha}(\zeta,\zeta)(z-\zeta)^\alpha(z_j-\zeta_j)\\
& +\sum_{j=1}^n\sum_{|\alpha|\leq k}\frac1{\alpha!(|\alpha|+1)} \diffp{^{|\alpha|+1}f}{\zeta^\alpha\partial\zeta_j}(z-\zeta)^\alpha(z_j-\zeta_j).
\end{align*}
Using (\ref{eq3}) we compute
\begin{align*}
 &\sum_{j=1}^n\sum_{|\alpha|\leq k}\frac1{\alpha!}\diffp{^{|\alpha|} \tilde b_j}{z^\alpha}(\zeta,\zeta)(z-\zeta)^\alpha(z_j-\zeta_j)\\
 &=\sum_{j=1}^n\sum_{|\alpha|\leq k}\frac{\alpha_j+1}{\alpha_1!\ldots(\alpha_j+1)!\ldots\alpha_n!}\diffp{^{|\alpha|} \tilde b_j}{z^{\alpha}}(\zeta,\zeta)(z_1-\zeta_1)^{\alpha_1}\ldots(z_j-\zeta_j)^{\alpha_j+1}\ldots (z_n-\zeta_n)^{\alpha_n}\\
 &=\sum_{j=1}^n\sum_{\over{|\alpha|\leq k+1}{\alpha_j\geq 1}}\frac{\alpha_j}{\alpha!}\diffp{^{|\alpha|-1} \tilde b_j}{z^{(\alpha_1,\ldots,\alpha_j-1,\ldots,\alpha_n)}}(\zeta,\zeta)(z-\zeta)^\alpha\\
 &=\sum_{{0<|\alpha|\leq k+1}}\sum_{j=1}^n\frac{\alpha_j}{\alpha!}\diffp{^{|\alpha|-1} \tilde b_j}{z^{(\alpha_1,\ldots,\alpha_j-1,\ldots,\alpha_n)}}(\zeta,\zeta)(z-\zeta)^\alpha\\
 &= \sum_{0<|\alpha|\leq k+1}\frac1{\alpha!} \diffp{^{|\alpha|}f}{\zeta^\alpha}(\zeta)(z-\zeta)^\alpha.
\end{align*}
We also compute
\begin{align*}
 &\sum_{j=1}^n \sum_{|\alpha|\leq k} \frac1{\alpha !(|\alpha|+1)}\diffp{^{|\alpha|+1}f}{\zeta^\alpha\partial\zeta_j}(\zeta)(z-\zeta)^\alpha(z_j-\zeta_j)\\
 &=\sum_{j=1}^n \sum_{|\alpha|\leq k} \frac{\alpha_j+1} {(|\alpha|+1)}\diffp{^{|\alpha|+1}f}{\zeta^{(\alpha_1,\ldots,\alpha_j+1,\ldots,\alpha_n)}}(\zeta)\frac{(z_1-\zeta_1)^{\alpha_1}\ldots(z_j-\zeta_j)^{\alpha_j+1}\ldots(z_n-\zeta_n)^{\alpha_n}}{\alpha_1 !\ldots (\alpha_j+1) !\ldots \alpha_n !}\\
 &=\sum_{j=1}^n \sum_{0<|\alpha|\leq k+1} \frac{\alpha_j}{\alpha ! |\alpha|}\diffp{^{|\alpha|}f}{\zeta^\alpha}(\zeta)(z-\zeta)^\alpha\\
 &=\sum_{0<|\alpha|\leq k+1} \frac{1}{\alpha !}\diffp{^{|\alpha|}f}{\zeta^\alpha}(\zeta)(z-\zeta)^\alpha
\end{align*}
from which we deduce the equality $\sum_{j=1}^n b_j(\zeta,z)(z_j-\zeta_j)=f(z)-f(\zeta)$.\qed
\section{Estimates of the extension operator}\label{section2}
We prove in this section that the previously defined extension operator satisfies the conclusion of Theorem \ref{th0}.
\subsection{More estimates}
In order to prove the $BMO$-estimates of Theorem \ref{th0} we apply the following classical
lemma:
\begin{lemma}\label{lemme-BMO}
 Let $h$ be a function of class $C^1$ on $D$. If there exists $C>0$ such that ${\rm d}h(\zeta)\leq
C|\rho(\zeta)|^{-1}$ then $h$ belongs to ${BMO}(D)$ and $\|h\|_{BMO(D)}\leq C$.
\end{lemma}
Thus we need to estimate the extension operator and its first derivatives. As usually with Berndtsson-Andersson kernel, problems occur when $\zeta$ and $z$ are close from each other and close from the boundary. Thus, it suffices to consider a point $z$ near $bD$, i.e. $z\in D\setminus\overline{D_{-\varepsilon_0}}$, and to integrate only for $\zeta$ in $\kb{\frac{\varepsilon_0}{2C}}z$. Moreover, the only interesting case for us is when $z$ is near a singularity of $X$.\\
We use the same notation as in Subsection \ref{Local parametrization} and we assume that $z=0$ belongs to $\sing(X)\cap bD$, and is a singularity of order $k$. We set $\omega=e^{\frac{2i\pi}k}$ and we assume that we are given a $c$-covering $(\kb{c|\rho(z_j)|}{z_j})_{j\in\nn}$ of a neighborhood $\cu(0)$ of $0$ where Propositions \ref{prop2.3}, \ref{prop2.4} and \ref{propII.0.1} hold true.
\begin{lemma}\label{lemme3.1}
 Let $p$ be any of the points among the $z_j$'s, $j\in\nn$. Then for all $\zeta'\in \pi_p(\kbz p)$, we have uniformly with respect to $\zeta'$ and $p$ :
 \begin{align*}
  \left|\frac{\prod_{j\in J_p}\left(\zeta'_n-\zeta'_1\varphi_p\left({\zeta'_1}^{\frac1k}\omega^j,\zeta'_2,\ldots,\zeta'_{n-1}\right)\right)}{f(\pi^{-1}_p(\zeta')) } \left(\pi_p^{-1}\right)^*b(\zeta',z)\right|&\leqs |\rho(p)|^{\frac{\#J_p-1}2}\sum_{\alpha=0}^k \left(\frac{\delta(p,z)}{|\rho(p)|}\right)^{\frac\alpha2},\\
  \left|{\rm d}_z\hskip -1pt\frac{\prod_{j\in J_p} \hskip -2pt\left(\zeta'_n-\zeta'_1\varphi_p\left(\hskip -1pt{\zeta'_1}^{\frac1k}\omega^j,\zeta'_2,\ldots,\zeta'_{n-1}\hskip -1pt\right)\hskip -1pt\right)}{f(\pi^{-1}_p(\zeta')) } \hskip -1pt\left(\hskip -1pt\pi_p^{-1}\hskip -1pt\right)^*\hskip -2pt b(\zeta',z)\right|&\leqs |\rho(p)|^{\frac{\#J_p-1}2-1}\hskip -1pt\sum_{\alpha=0}^k \hskip -2pt\left(\frac{\delta(p,z)}{|\rho(p)|}\right)^{\frac\alpha2}.
 \end{align*}
\end{lemma}
\pr We only prove the first inequality, the second can be proved in the same way.
Since for all $\zeta'\in \pi_p(\kbz p)$ we have $\delta(\pi_p^{-1}(\zeta'),z)\leqs \delta(\pi_p^{-1}(\zeta'),p)+\delta(p,z)\leqs |\rho(p)|+\delta(p,z) $, it suffices to prove the inequality with $\delta(\pi_p^{-1}(\zeta'),z)$ instead of $\delta(p,z)$.

We have to distinguish two cases. 
\begin{trivlist}
 \item[First case :] When $d(p,\sing(X))\leq 10C|c\rho(p)|^{\frac12}$, according do its definition, $J_p=\{0,\ldots, k-1\}$, and so
 we have to prove that 
  \begin{align}
\left|\left(\pi_p^{-1}\right)^*b(\zeta',z)\right|&\leq |\rho(p)|^{\frac {k-1}2}\sum_{\alpha=0}^k \left(\frac{\delta(\pi_p^{-1}(\zeta'),z)}{|\rho(p)|}\right)^{\frac\alpha2}.   \label{eq1bis}
  \end{align}
We prove that for $K>0$ and all $\zeta\in B(p,K|\rho(p)|^{\frac12})$, we have $|f(\zeta)|\leqs |\rho(p)|^{\frac{k}{2}}$, uniformly with respect to $p$ and $\zeta$.
We have $d(\zeta,\sing(X))\leq d(p,\sing(X))+|\zeta-p|\leqs |\rho(p)|^{\frac12}$. Since $d(\zeta,\sing(X))=|\zeta_1|+|\zeta_n|$, setting $\zeta'=\pi_p(\zeta)$, it comes $|\zeta'_1|\leqs|\rho(p)|^{\frac12}$ and $|\zeta'_n|\leqs |\rho(p)|^{\frac12}$.
Now, if we denote by $(\zeta'_1)^{\frac1k}$ a complex number such that $\left((\zeta'_1)^{\frac1k}\right)^k=\zeta'_1$, for any $j$ we have
\begin{align*}
 \left|\zeta'_n-\zeta'_1\varphi_p({\zeta'_1}^\frac1k\omega^j,\zeta'_2,\ldots,\zeta'_{n-1})\right|\leqs |\rho(p)|^{\frac12}.
\end{align*}
This readily implies that $|f(\zeta)|\leqs |\rho(p)|^{\frac k2}$.\\
It then comes for all multi-index $\alpha$ and all $\zeta \in  B(p,K|\rho(p)|^{\frac12})$ the inequality $\left|\diffp{^{|\alpha|} f}{\zeta^\alpha}(\zeta)\right| \leqs |\rho(p)|^{\frac {k-|\alpha|}2}$. 
Now Inequality (\ref{eq1bis}) then immediately follows from lemma \ref{lemme Hefer}.

 \item[Second case :] When $d(p,\sing(X))\geq 10C|c\rho(p)|^{\frac12}$, we put $f^*_p=f\circ\pi^{-1}_p$ and $b^*=\left(\pi_p^{-1}\right)^* b=\sum_{j=1}^nb^*_jd\zeta'_j$ where $\zeta'=\pi_p(\zeta)$. For all $l$ we have
 \begin{align*}
  b^*_{l}(\zeta',z)&=\sum_{|\alpha|\leq k}\frac1{\alpha!(|\alpha|+1)} \diffp{^{|\alpha|+1}f^*_p}{{\zeta'}^\alpha\partial\zeta'_l}(\zeta')(\pi_p(z)-\zeta')^\alpha +O(|z-\pi^{-1}_p(\zeta')|^{k+1}).
 \end{align*}
Again, we have to control the derivatives of $f$.
Proposition \ref{prop2.2} gives 
\begin{align}
f^*_p(\zeta')&=u_p(\pi_p^{-1}(\zeta')) \prod_{j=0}^{k-1}\left(\zeta'_n-\zeta'_1\varphi_p({\zeta'_1}^{\frac1k}\omega^j,\zeta'_2,\ldots,\zeta'_{n-1})\right)\label{eq1ter}
\end{align}
where $\xi\mapsto \xi^{\frac1k}$ is a holomorphic $k$-root defined on $D(\pi_{p,1}(p),C|c\rho(p)|^{\frac12})$.

For simplicity sake, we will assume that $u_p$ is in fact constant because, since $u_p$ and its derivatives are uniformly bounded, $u_p$ won't play any role.\\
We deduce from equality (\ref{eq1ter}) that for all $l$, $\left|\diffp{^{|\alpha|+1}f_p^*}{{\zeta'}^\alpha\partial\zeta'_l}(\zeta')\right|$ is smaller than a sum over all sets $F\subset\{0,\ldots, k-1\}$ such that $\#F\geq k-|\alpha|-1$  of the following terms 
\begin{align*}
 \prod_{j\in F} \left|\zeta'_n-\zeta'_1\varphi_p({\zeta'_1}^{\frac1k}\omega^j,\zeta'_2,\ldots,\zeta'_{n-1})\right|
 \prod_{j\in\{0,\ldots, k-1\}\setminus F } \diffp{^{|\beta_j|} \left(\zeta'_n-\zeta'_1\varphi_p({\zeta'_1}^{\frac1k}\omega^j,\zeta'_2,\ldots,\zeta'_{n-1})\right)}{{\zeta'}^{\beta_j}}.
\end{align*}
where the $\beta_j$'s, $j\in\{0,\ldots, k-1\}\setminus F$, are  multi-index such that $\sum_j \beta_j=\alpha +(0,\ldots,1,\ldots, 0)$, the $1$ being at $l$-th position.\\
For all $\zeta'\in\pi_p(\kbz p),$ from Proposition \ref{prop2.4}, $\zeta'_1$ belongs to  $D(\pi_{p,1}(p),C|c\rho(p)|^{\frac12})$ and since $|\pi_{p,1}(p)|\geq 2C|c\rho(p)|^{\frac12}$, we have $|\zeta'_1|\geqs |\rho(p)|^{\frac12}$. Since from Lemma \ref{lemme2.1}, the derivatives of $\varphi_p$ are bounded, we get
\begin{align*}
\left|\diffp{^{|\beta_j|} \left(\zeta'_n-\zeta'_1\varphi_p({\zeta'_1}^{\frac1k}\omega^j,\zeta'_2,\ldots,\zeta'_{n-1})\right)}{{\zeta'}^{\beta_j}}\right|&\leqs |\rho(p)|^{-\frac{|\beta_j|-1}2}
\end{align*}
and then, denoting by $A^c$ the complement of the set $A$ in $\{0,\ldots,k-1\}$, we have :
\begin{align*}
&\left|\frac{\prod_{j\in J_p} \left(\zeta'_n-\zeta'_1\varphi_p({\zeta'_1}^{\frac1k}\omega^j,\zeta'_2,\ldots,\zeta'_{n-1})\right)}{f^*_p(\zeta')}
\diffp{^{|\alpha|+1}f^*_p}{{\zeta'}^\alpha\partial\zeta'_l}(\zeta')\right|\\
&\leqs 
\sum_{\over{F\subset\{0,\ldots, k-1\}}{\#F\leq k-|\alpha|-1}}
\frac{\prod_{j\in F}\left(\zeta'_n-\zeta'_1\varphi_p({\zeta'_1}^{\frac1k}\omega^j,\zeta'_2,\ldots,\zeta'_{n-1})\right)}{\prod_{J_p^c}\left(\zeta'_n-\zeta'_1\varphi_p({\zeta'_1}^{\frac1k}\omega^j,\zeta'_2,\ldots,\zeta'_{n-1})\right)}|\rho(p)|^{-\frac{1+|\alpha|-k+\#F}2} \delta(\pi^{-1}_p(\zeta'),z)^{\frac{|\alpha|}{2}}\\
&\leqs\sum_{\over{F\subset\{0,\ldots, k-1\}}{\#F\leq k-|\alpha|-1}}
\frac{\prod_{j\in F\cap J_p}\left(\zeta'_n-\zeta'_1\varphi_p({\zeta'_1}^{\frac1k}\omega^j,\zeta'_2,\ldots,\zeta'_{n-1})\right)}{\prod_{J_p^c\cap F^c}\left(\zeta'_n-\zeta'_1\varphi_p({\zeta'_1}^{\frac1k}\omega^j,\zeta'_2,\ldots,\zeta'_{n-1})\right)}|\rho(p)|^{-\frac{1+|\alpha|-k+\#F}2} \delta(\pi^{-1}_p(\zeta'),z)^{\frac{|\alpha|}{2}}.
\end{align*}
The definition of $J_p$ implies that $\left|\zeta'_n-\zeta'_1\varphi_p({\zeta'_1}^{\frac1k}\omega^j,\zeta'_2,\ldots,\zeta'_{n-1})\right|\geqs|\rho(p)|^{\frac12}$ for all $\zeta'\in\pi_p(\kbz p)$ and all $j\notin J_p$. From Corollary \ref{cor2.5} comes $\left|\zeta'_n-\zeta'_1\varphi_p({\zeta'_1}^{\frac1k}\omega^j,\zeta'_2,\ldots,\zeta'_{n-1})\right|\leqs|\rho(p)|^{\frac12}$ for all $\zeta'\in\pi_p(\kbz p)$ and all $j\in J_p$. Therefore
\begin{align*}
&\left|\frac{\prod_{j\in J_p} \left(\zeta'_n-\zeta'_1\varphi_p({\zeta'_1}^{\frac1k}\omega^j,\zeta'_2,\ldots,\zeta'_{n-1})\right)}{f^*_p(\zeta')}
\diffp{^{|\alpha|+1}f_p^*}{{\zeta'}^\alpha\partial\zeta'_l}(\zeta')\right|\\
&\leqs\sum_{\over{F\subset\{0,\ldots, k-1\}}{\#F\leq k-|\alpha|-1}}
|\rho(p)|^{\frac{\#(F\cap J_p)-\#(J_p^c\cap F^c)-1-|\alpha|+k-\#F}2} \delta(\pi^{-1}_p(\zeta'),z)^{\frac{|\alpha|}{2}}.
\end{align*}

Since $\#(F\cap J_p)+\#F^c=\#J_p+\#(J_p^c\cap F^c)$, this leads to
\begin{align*}
\left|\frac{\prod_{j\in J_p} \left(\zeta'_n-\zeta'_1\varphi_p({\zeta'_1}^{\frac1k}\omega^j,\zeta'_2,\ldots,\zeta'_{n-1})\right)}{f_p^*(\zeta')}
\diffp{^{|\alpha|+1}f_p^*}{{\zeta'}^\alpha\partial\zeta'_l}(\zeta')\right|&\leqs |\rho(p)|^{\frac{\#J_p-1}2} \hskip -1pt\left(\frac{\delta(\pi^{-1}_p(\zeta'),z)}{|\rho(p)|}\right)^{\frac{|\alpha|}{2}},
\end{align*}
which was to be shown.\qed
\end{trivlist}
\begin{lemma}
 There exists a neighborhood $\cu(bD)$ of $bD$ such that for all $p\in\cu(bD)\cap D$, all $\zeta\in\kbz p \cap D$ and all $z\in D$, we have
 $$|1+G(\zeta,z)|\geqs \frac1{|\rho(\zeta)|}(|\rho(z)|+|\rho(\zeta)|+\delta(\zeta,z)).$$
\end{lemma}
\pr For all $p$, $\zeta\in\kbz p$ we have $\delta(p,\zeta)\leq c|\rho(p)|$ and $|\rho(\zeta)|\geq \frac12 |\rho(p)|$. Choosing $c$  such that $c<\frac14$, it comes :
\begin{align*}
 |\rho(\zeta)|+\delta(z,\zeta)
 &\geq \frac12|\rho(p)|+\frac1{C} \delta(z,p)-\delta(\zeta,p)\\
 &\geq \frac12 |\rho(p)|+\frac1{C} \delta(p,z) -c|\rho(p)|\\
 &\geqs |\rho(p)|+\delta(p,z)
\end{align*}
and therefor it suffices to prove that $|\rho(\zeta)+g^{(0)}(\zeta,z)|\geqs |\rho(z)|+|\rho(\zeta)|+\delta(\zeta,z)$. Moreover, it suffices to consider the case $|\zeta-z|<\varepsilon$ for some arbitrary small $\varepsilon>0$.
From Subsection \ref{subsection2.5} we have
\begin{align*}
g^{(0)}&=\frac{(F+|F|^2(\overline v -\sup_{\cu(\overline{D})^2}|v|)} {2|1+F(v-\sup_{\cu(\overline{D})^2}|v|)|^2}
\end{align*}
and so
\begin{align*}
\re\ g^{(0)}&=\frac{(\re F+|F|^2 \re(\overline v -\sup_{\cu(\overline{D})^2}|v|)} {2|1+F(v-\sup_{\cu(\overline{D})^2}|v|)|^2}.
\end{align*}
If $|\zeta-z|$ is sufficient small, $|1+F(v-\sup_{\cu(\overline{D})^2}|v|)|^2\geq \frac12$ and so
\begin{align*}
\rho(\zeta)+\re\ g^{(0)}(\zeta,z)&\leq \rho(\zeta)+\re\ F(\zeta,z)\\
&\leq \frac{\rho(z)+\rho(\zeta)}2-\beta|\zeta-z|^2.
\end{align*}
This implies that 
\begin{align}
|\rho(\zeta)+\re\ g^{(0)}(\zeta,z)|&\geq \frac{|\rho(z)|+|\rho(\zeta)|}2+\beta|\zeta-z|^2.\label{eq24}
\end{align}
We also have $\delta(\zeta,z)\eqs |\zeta-z|^2+|\langle \eta_\zeta,z-\zeta\rangle|$. So if $|\langle \eta_\zeta,z-\zeta\rangle|\leq \gamma|\zeta-z|^2$ for some arbitrarily big $\gamma>0$, (\ref{eq24}) gives $|\rho(\zeta)+\re\ g^{(0)}(\zeta,z)|\geqs |\rho(z)|+|\rho(\zeta)|+\delta(\zeta,z)$.\\
If $|\langle \eta_\zeta,z-\zeta\rangle|\geq \gamma|\zeta-z|^2$ 
\begin{align*}
 |g^{(0)}(\zeta,z)|
 &\geqs |\re\ g^{(0)}(\zeta,z)|+ |\im\ g^{(0)}(\zeta,z)|\\
 &\geq |F(\zeta,z)|-O(|\zeta-z|^2)\\
 &\geqs |\langle \eta_\zeta,z-\zeta\rangle| -O(|\zeta-z|^2)\\
 &\geqs \delta(\zeta,z)
\end{align*}
 provided $\gamma$ is big enough. With (\ref{eq24}), we are done in this case too.\qed
%

We denote by $\psi_p$ the smooth cutoff function of the partition of unity related to $\kbz p$. The following corollary then immediately follows :\begin{corollary}\label{corollaire3.4}
 There exists a neighborhood $\cu(bD)$ of $bD$  such that for all $p\in\cu(bD)\cap D$, all $\zeta\in\kbz p \cap D$ and all $z\in D$ :
 \begin{align*}
\left|\overline\partial \diffp{^{\#J_p}}{{\zeta'_n}^{\#J_p}}\left(
{\pi_p^{-1}}^* \left(\psi_p P^{N,n-1}\right)(\zeta',z)
\right)\right|&\leqs \left(\frac{|\rho(p)|}{|\rho(p)|+\delta(p,z)}\right)^N \frac1{|\rho(p)|^{n+\frac12+\frac{\#J_p}2}},\\
\left|d_z\overline\partial \diffp{^{\#J_p}}{{\zeta'_n}^{\#J_p}}\left(
{\pi_p^{-1}}^* \left(\psi_p P^{N,n-1}\right)(\zeta',z)
\right)\right|&\leqs \left(\frac{|\rho(p)|}{|\rho(p)|+\delta(p,z)}\right)^N \frac1{|\rho(p)|^{n+\frac32+\frac{\#J_p}2}}.
\end{align*}
\end{corollary}
\subsection{BMO-extension}
Let $\tilde h$ be a smooth extension of $h$ as in the hypothesis of Theorem \ref{th0}.
We set $\gamma_\infty=\sup_{\over{\zeta\in D}{|\alpha|\leq k}}
\left|\diffp{^{\alpha}\tilde h}{\overline{\epsilon_1}^{\alpha_1}\ldots \partial \overline{\epsilon_n}^{\alpha_n}}(\zeta)\right||\rho(\zeta)|^{\alpha_1+\frac{\alpha_2+\ldots+\alpha_n}2}$.
In order to prove Theorem \ref{th0} when $q=+\infty$, we have to prove that $E_N (g)$ is in $BMO(D)$ and $\|E_N (g)\|_{BMO(D)}\leqs
\gamma_\infty$.\\
We keep the notations of the previous section. If $p$ is any point among the points $z_j$ of the covering, we get from Lemma \ref{lemme3.1} and Corollary \ref{corollaire3.4} for $\zeta'\in \pi_p^{-1}(\kbz p)$ and $z\in D$ :
\begin{align*}
& \left|d_z
 \left(\frac{\prod_{j\in J_p}\left(\zeta'_n-\zeta'_1\varphi_p\left({\zeta'_1}^{\frac1k}\omega^j,\zeta'_2,\ldots,\zeta'_{n-1}\right)\right)}{f(\pi^{-1}_p(\zeta')) } \overline\partial \diffp{^{\#J_p}}{{\zeta'_n}^{\#J_p}}\left(
{\pi_p^{-1}}^* \left(\psi_p h b\wedge P^{N,n-1}\right)(\zeta',z)
\right)\right)\right|\\
&\leqs \gamma_{\infty}\left(\frac{|\rho(p)|}{|\rho(p)|+\delta(p,z)}\right)^{N'} \frac1{|\rho(p)|^{n+2}}.
\end{align*}
for some integer $N'\geq 2 $. Therefor, when we integrate over $\kbz p$, we get
\begin{align*}
& \Bigg| d_z\int_{\zeta'\in\pi_p(\kbz p)}
 \frac{\prod_{j\in J_p}\left(\zeta'_n-\zeta'_1\varphi_p\left({\zeta'_1}^{\frac1k}\omega^j,\zeta'_2,\ldots,\zeta'_{n-1}\right)\right)}{f(\pi^{-1}_p(\zeta')) }\\
 &\overline\partial
 \diffp{^{\#J_p}}{{\zeta'_n}^{\#J_p}}
 \left({
	\pi_p^{-1}}^* 
	 \left(
	      \psi_p h b\wedge P^{N,n-1}
	 \right)
	 (\zeta',z)
\right)\Bigg|\\
&\leqs \gamma_{\infty} \frac{|\rho(p)|^{N'-1}}{\left(|\rho(p)|+\delta(p,z)\right)^{N'}}
\end{align*}
Now, we renumber the covering $\kbc$ of $\cu(0)$ as in Lemma \ref{nombre_boules}. In order to apply Lemma \ref{lemme-BMO}, we  show that 
$\sum_{j=j_0}^\infty\sum_{i=0}^{i_0(j)}\sum_{m=1}^{m_{i,j}}
\frac{|\rho(z_{m}^{i,j})|^{N'-1}}{\left((i+1)|\rho(z^{i,j}_m)|+|\rho(z)|\right)^{N'}}
$ is uniformly bounded by $\frac1{|\rho(z)|}$. We have:
\begin{eqnarray*}
\lefteqn{\sum_{j=j_0}^\infty\sum_{i=0}^{i_0(j)}\sum_{m=1}^{m_{i,j}}\frac{|\rho(z_{m}^{i,j})|^{N'-1}}{\left((i+1)|\rho(z^{i,j}_m)|+|\rho(z)|\right)^{N'}}}\\
&\leq&
 \sum_{j=j_0}^\infty\sum_{i=0}^{i_0(j)}\sum_{m=1}^{m_{i,j}}
\left(\frac{(1-c\kappa)^j}{(i+1)(1-c\kappa)^j+1}\right)^{N'-1} \cdot\frac{1}{((i+1)(1-c\kappa)^j+1)|\rho(z)|}\\
&\leq&\frac{1}{|\rho(z)|}\left(\sum_{j=0}^\infty \sum_{i=0}^\infty
\frac{(1-c\kappa)^j}{(i+1)^{N'-2-n}} +
\sum_{j=j_0}^{-1} \sum_{i=0}^\infty \frac1{(i+1)^{N'-n}(1-c\kappa)^{j}}\right)\\
&\leqs& \frac1{|\rho(z)|}.
\end{eqnarray*}
So $E_N(g)$ belongs to ${BMO}(D)$ and $\|E_N(g)\|_{BMO(D)}$, up to a multiplicative uniform constant, is lower than 
$\sup_{\over{\zeta\in D}{|\alpha|\leq k}}
\left|\diffp{^{\alpha}\tilde h}{\overline{\epsilon_1}^{\alpha_1}\ldots \partial \overline{\epsilon_1}^{\alpha_n}}(\zeta)\right||\rho(\zeta)|^{\alpha_1+\frac{\alpha_2+\ldots+\alpha_n}2}$.\qed
\subsection{$L^q$ extension}
The $L^q$-estimates of Theorem \ref{th0} are left to be shown. For $q\in ]1,+\infty[$ we will apply the following lemma (see
\cite{Pol}):
\begin{lemma}\mlabel{Pol}
 Suppose the kernel $k(\zeta,z)$ is defined on $D\times D$ and the operator $K$
is defined by $Kf(z)=\int_{\zeta\in D}k(\zeta,z)f(\zeta)d\lambda(\zeta)$. If for
every $\varepsilon\in ]0,1[$, there exists a constant $c_\varepsilon$ such that
\begin{eqnarray*}
 \int_{\zeta\in D} |\rho(\zeta)|^{-\varepsilon}|k(\zeta,z)|d\lambda(\zeta)&\leq&
c_\varepsilon |\rho(z)|^{-\varepsilon},\quad \forall z\in D,
\end{eqnarray*}
and
\begin{eqnarray*}
\int_{z\in D} |\rho(z)|^{-\varepsilon}|k(\zeta,z)|d\lambda(z)&\leq& c_\varepsilon
|\rho(\zeta)|^{-\varepsilon},\quad \forall \zeta\in D,
\end{eqnarray*}
then for all $q\in ]1,+\infty[$, there exists $c_q>0$ such that
$\|Kf\|_{L^q(D)}\leq \|f\|_{L^q(D)}$.
\end{lemma}
{\it Proof of Theorem \ref{th0} for $q\in]1,+\infty[$~:}
Applying Lemma \ref{lemme3.1} and \ref{Pol} and Corollary \ref{corollaire3.4}, it
suffices to prove that for all $\varepsilon\in ]0,1[$ there exists
$c_\varepsilon>0$ such that
\begin{align}
{\int_{\zeta\in D}
\frac {|\rho(\zeta)|^{N'-\varepsilon}}
{\left(|\rho(\zeta)|+|\rho(z)|+\delta(\zeta,z)\right)^{N'+n+1}}
d\lambda(\zeta)}&\leq c_{\varepsilon}|\rho(z)|^{-\varepsilon},\ \forall z\in D,\mlabel{eq5}\\
{\int_{z\in D}
\frac {|\rho(\zeta)|^{N'}|\rho(z)|^{-\varepsilon}}
{\left(|\rho(\zeta)|+|\rho(z)|+\delta(\zeta,z)\right)^{N'+n+1}}
d\lambda(z)}&\leq c_{\varepsilon}|\rho(\zeta)|^{-\varepsilon},\ \forall \zeta \in D.\mlabel{eq6}
\end{align}
The inequality (\ref{eq5}) can be shown as in the proof of Theorem \ref{th0} for $q=\infty$.\\
In order to prove that the inequality (\ref{eq6}) holds true, we cover $D$ with the \ko balls $\kbz\zeta$ and $\left({\cal P}_{2^{j+1}c|\rho(\zeta)|}(\zeta)\setminus {\cal
P}_{2^jc|\rho(\zeta)|}(\zeta)\right)$, $j\in\nn$.\\
For $z\in{\cal P}_{c|\rho(\zeta)|}(\zeta )$, $|\rho(z)|\eqs |\rho(\zeta)|$ and thus
\begin{align}
\int_{z\in{\cal P}_{c|\rho(\zeta)|}(\zeta )} \frac {|\rho(\zeta)|^{N'}|\rho(z)|^{-\varepsilon}}
{\left(|\rho(\zeta)|+|\rho(z)|+\delta(\zeta,z)\right)^{N'+n+1}}d\lambda(z)
&\leqs|\rho(\zeta)|^{-\varepsilon}.\mlabel{eq20}
\end{align}
When we integrate on ${\cal P}_{2^{j+1}c|\rho(\zeta)|}(\zeta)\setminus {\cal
P}_{2^jc|\rho(\zeta)|}(\zeta)$ we get
\begin{align}
&\nonumber \int_{{\cal P}_{2^{j+1}c |\rho(\zeta)|}(\zeta)\setminus {\cal
P}_{2^jc |\rho(\zeta)|}(\zeta)} 
 \frac {|\rho(\zeta)|^{N'}|\rho(z)|^{-\varepsilon}}
{\left(|\rho(\zeta)|+|\rho(z)|+\delta(\zeta,z)\right)^{N'+n+1}}d\lambda(z)\\
\nonumber &\leqs
\int_{\over{|x_1|,|y_1|\leq 2^{j+1}c|\rho(\zeta)|}{|x_2|,|y_2|,\ldots,|x_n|,|y_n|\leq \sqrt{2^{j+1}c|\rho(\zeta)|}}}
 \frac{|\rho(\zeta)|^{N'}x_1^{-\varepsilon}}
{\left(|\rho(\zeta)|+2^jc|\rho(\zeta)|\right)^{N'+n+1}}d\lambda(z)\\
\nonumber&\leqs
(2^{j+1}c|\rho(\zeta)|)^{-\varepsilon+n+1}
\frac{|\rho(\zeta)|^{N'}}
{\left(|\rho(\zeta)|+2^jc|\rho(\zeta)|\right)^{N'+n+1}}\\
&\leqs |\rho(\zeta)|^{-\varepsilon} 2^{-j(N'+\varepsilon)}.\mlabel{eq21}
\end{align}
Summing  (\ref{eq20}) and (\ref{eq21}) for all non-negative integer $j$ we prove inequality (\ref{eq20}).  Theorem \ref{th0} is therefore proved
for $q\in]1,+\infty[$.\qed\\
{\it Proof of Theorem \ref{th0} for $q=1$~:} We prove directly that $E_N(g)$ belongs to $L^1(D)$. 
Lemma \ref{lemme3.1} and Corollary \ref{corollaire3.4} yield
\begin{align*}
{\int_D\hskip -2pt |E_N g(z)|d\lambda(z)}\hskip -1pt &\leqs \hskip -1pt
\sum_{j=0}^\infty\sum_{|\alpha|\leq \#J_{z_j}+1} 
\int_{\kbz {z_j}}\hskip -3pt  {|\rho(z_j)|^{\alpha_1+\frac{\alpha_2+\ldots+\alpha_n}2}}\hskip -3pt
\left|\hskip -1pt \diffp{^{|\alpha|}\tilde h} {\overline{\epsilon_1(z_j)}^{\alpha_1}\ldots \partial \overline{\epsilon_n(z_j)}^{\alpha_n}}(\zeta)\hskip -1pt\right|\\
& \cdot \left(\int_D
\frac{|\rho(\zeta)|^{N'}}{\left(|\rho(\zeta)|+|\rho(z)|+\delta(\zeta,z)\right)^{N'+n+1}}
{d\lambda(z)}\right)d\lambda(\zeta).
\end{align*}
We may show that $\int_D
\frac{|\rho(\zeta)|^{N'}}{\left(|\rho(\zeta)|+|\rho(z)|+\delta(\zeta,z)\right)^{N'+n+1}}
{d\lambda(z)}$ is bounded exactly as (\ref{eq6}) and we don't repeat it here. We then get
\begin{align*}
&{\int_D|E_Ng(z)|d\lambda(z)}\\
&\leqs\sum_{j=0}^\infty\sum_{ |\alpha|\leq \#J_{z_j}+1} \int_{\kbz {z_j}} 
{|\rho(z_j)|^{\alpha_1+\frac{\alpha_2+\ldots+\alpha_n}2}}
\left| \diffp{^{|\alpha|}\tilde h} {\overline{\epsilon_1(z_j)}^{\alpha_1}\ldots \partial \overline{\epsilon_n(z_j)}^{\alpha_n}}(\zeta)\right|d\lambda(\zeta)\\
&\leqs
\sum_{0\leq |\alpha|\leq k+1}
\left\| \zeta \mapsto
{|\rho(\zeta)|^{\alpha_1+\frac{\alpha_2+\ldots+\alpha_n}2}}
\diffp{^{\alpha}\tilde h}{\overline{\epsilon_1}^{\alpha_1}\ldots \partial \overline{\epsilon_n}^{\alpha_n}}(\zeta)
 \right\|_{L^1(D)}.
\end{align*}
\qed
\section{Smooth $L^q$ extension}\label{section3}
\subsection{Case $q=+\infty$}
In this subsection, we aim at proving Theorem \ref{th1}. We choose $c>0$ sufficiently small so that $\kb{c|\rho(z)|}z$ is included in $D_{\frac{\rho(z)}2}$ for every $z$ sufficiently close to $bD$, and we prove the following lemma.
\begin{lemma}\label{lemme4.1}
  Let $h\in\co(X\cap D)$ be such that for all $z_0$ sufficiently close to $bD$, there exists an extension $h_0\in\co(\kb{c|\rho(z_0)|}{z_0})$ of $h$ satisfying $\sup_{\kb{c|\rho(z_0)|}{z_0}}|h_0|\leqs 1$, uniformly with respect to $z_0$.\\
 Then there exists an extension $\tilde h\in C^\infty(D)$ of $h$ which satisfies hypothesis (\ref{th0i}-\ref{th0iii}) of Theorem \ref{th0} for $q=+\infty$.
\end{lemma}
\pr We choose a $\frac c2$-covering $\left(\kb{\frac c2|\rho(z_j)|}{z_j}\right)_{j\in\nn^*}$ of $D\setminus \overline{D_{-\varepsilon_0}}$ where $\varepsilon_0>0$ is  sufficiently small. Let $\epsilon_1(z_j)=\eta_{z_j},\ldots,\epsilon_n(z_j)$ be an orthonormal basis of $\cc^n$. We denote by $\zeta^*=(\zeta^*_1,\ldots,\zeta^*_n)$ the coordinates system centered at $z_j$ and of basis $\epsilon_1(z_j),\ldots,\epsilon_n(z_j)$. We now choose a partition of unity $(\psi_j)_{j\in\nn}$ relative to the covering $\left(\kb{\frac c2|\rho(z_j)|}{z_j}\right)_{j\in\nn^*}\cup D_{-\varepsilon_0}$ of $D$ such that for all $j>0$, for all multi-indexes $\alpha$ and $\beta$, 
$\left|
\diffp{^{|\alpha|+|\beta|}\psi_j}{{\zeta^*}^\alpha\partial\overline{\zeta^*}^\beta}(\zeta)
\right|\leqs\frac{1}{|\rho(z_j)|^{\alpha_1+\beta_1+\frac{\sum_{l=2}^n\alpha_l+\beta_l}2}}.$\\
Now, we put $\tilde h=\sum_{j=0}^{+\infty}\psi_j h_j$ where $h_0$ is a holomorphic extension of $h$ given by Cartan B theorem and where, for $j>0$, $h_j$ is the holomorphic extension of $h$ on $\kbz{z_j}$ given by the hypothesis of the lemma.\\
Cauchy's inequalities imply that for all $j$, all $\zeta\in \kb{\frac c2|\rho(z_j)|}{z_j}$ and all multi-index $\alpha$, $\left|\diffp{^{|\alpha|}h_j}{{\zeta^*}^\alpha}(\zeta)\right|\leqs \frac1{|\rho(z_j)|^{\alpha_1+\frac{\alpha_2+\ldots+\alpha_n}2}}$ uniformly with respect to $\zeta$ and $j$.
Therefore, for all multi-indexes $\alpha$ and $\beta$, 
$\left|
\diffp{^{|\alpha|+|\beta|}\psi_jh_j}{{\zeta^*}^\alpha\partial\overline{\zeta^*}^\beta}(\zeta)
\right|\leqs\frac{1}{|\rho(z_j)|^{\alpha_1+\beta_1+\frac{\sum_{l=2}^n\alpha_l+\beta_l}2}}.$
It follows immediately that for all $N\in\nn^*$, $\rho^{N} \tilde h$ vanishes to order $N$ on $bD$.
Since $|\rho(\zeta)|\eqs|\rho(z_j)|$ on ${\cal P}_{\frac c2 |\rho(z_j)|}(z_j)$, for all multi-index $\alpha$ also comes the uniform boundedness of ${|\rho(\zeta)|^{\alpha_1+\frac{\sum_{l=2}^n\alpha_l}2}}\diffp{^{|\alpha|}\psi_jh_j}{\overline{\zeta^*}^\alpha}(\zeta)$ and so of
$\left|\diffp{^{\alpha}\tilde h}{\overline{\epsilon_1}^{\alpha_1}\ldots \partial \overline{\epsilon_n}^{\alpha_n}}\right||\rho|^{\alpha_1+\frac{\alpha_2+\ldots+\alpha_n}2}$.
Finally, since for all $j$ the functions $h_j$ are holomorphic, $\diffp{^{\alpha}\tilde h}{\overline{\epsilon_1}^{\alpha_1}\ldots \partial \overline{\epsilon_n}^{\alpha_n}}=0$ on $X\cap D$ for all
 multi-index $\alpha$.\qed

Therefor, in order to prove Theorem \ref{th1}, it suffices to extend locally and uniformly $h$. We will achieve this goal with divided differences. 
If $\varphi$ is a function defined on an open set $\cu$ of $\cc$ and if $t_1,\ldots, t_j$ are $j$ pairwise distinct points in $\cu$, we set
\begin{align*}
 \varphi[t_l]&=\varphi(t_l),\ l=1,\ldots, j,\\
 \varphi[t_1,\ldots, t_j]&=\frac{\varphi[t_2,\ldots,t_j]-\varphi[t_1,\ldots,t_{j-1}]}{t_j-t_1}.
\end{align*}
%

{\it Proof of Theorem \ref{th1} :} according to Lemma \ref{lemme4.1}, it suffices to extend $h$ holomorphically and boundedly on any Koranyi ball $\kb{c|\rho(p)|}{p}$ where $p\in D$ is close to $bD$. Without restriction, we assume that $\kb{c|\rho(p)|}{p}$ intersects $X$ and that $p$ belongs to a neighborhood $\cu(0)$ of the point $0$ which belongs to $bD\cap \sing(X)$. We adopt the notations of Subsection \ref{Local parametrization}. Let $k$ denotes the order of the singularity $0$.

If $d(p,\sing(X))\leq 10C|c\rho(p)|^{\frac12}$, we set $J_p=\{0,\ldots,k-1\}$ and for any $\xi$, we choose an arbitrary complex number $\xi^{\frac1k}$ such that $\xi=(\xi^{\frac1k})^k$.

If $d(p,\sing(X))\geq 10C|c\rho(p)|^{\frac12}$, we set as in Subsection \ref{Local parametrization} :
$$J_p=\left\{
j\in\{0,\ldots,k\}\tq \exists \zeta' \in\pi_p\big(\kbz p\big),\ \big|\zeta'_1\varphi_p({\zeta'_1}^{\frac1k}\omega^j,\zeta'_2,\ldots,\zeta'_{n-1})-\zeta'_n\big|\leq 2|c\rho(p)|^{\frac12}
\right\}$$
where $\xi\mapsto \xi^{\frac1k}$ is a $k$-root defined on $D(\pi_{p,1}(p),C|c\rho(p)|^{\frac12})$ and $\omega=e^{\frac{2i\pi}k}$.
We then define $h_p$ by setting $h_p=H_p\circ \pi_p^{-1}$ and
\begin{align*}
 H_{p}(\zeta')
 =&\sum_{j\in J_z} \left( \prod_{\over{l\in J_p}{l\neq j}}
\left( \frac{\zeta'_n-\zeta'_1\varphi_p({\zeta'_1}^{\frac1k}\omega^l,\zeta'_2,\ldots,\zeta'_{n-1})}
 {\zeta'_1\varphi_p({\zeta'_1}^{\frac1k}\omega^j,\zeta'_2,\ldots,\zeta'_{n-1})-\zeta'_1\varphi_p({\zeta'_1}^{\frac1k}\omega^l,\zeta'_2,\ldots,\zeta'_{n-1})}\right)\right.\\
 &\hphantom{ \sum_{j\in J_z} \left( \prod_{\over{l\in J_p}{l\neq j}}\right.}\cdot\left. 
 h  \Big(\pi^{-1}_p\big(\zeta'_1,\ldots,\zeta'_{n-1},\zeta'_1\varphi_p({\zeta'_1}^{\frac1k}\omega^l,\zeta'_2,\ldots,\zeta'_{n-1})\big)\Big)
 \vphantom{\prod_{\over{l\in J_p}{l\neq j}}
\left( \frac{\zeta'_n-\zeta'_1\varphi_p({\zeta'_1}^{\frac1k}\omega^l,\zeta'_2,\ldots,\zeta'_{n-1})}
 {\zeta'_1\varphi_p({\zeta'_1}^{\frac1k}\omega^j,\zeta'_2,\ldots,\zeta'_{n-1})-\zeta'_1\varphi_p({\zeta'_1}^{\frac1k}\omega^l,\zeta'_2,\ldots,\zeta'_{n-1})}\right)}\right).
\end{align*}
Thus $H_p(\zeta')$ is well defined provided that $\pi^{-1}_p\big(\zeta'_1,\ldots,\zeta'_{n-1},\zeta'_1\varphi_p({\zeta'_1}^{\frac1k}\omega^l,\zeta'_2,\ldots,\zeta'_{n-1}))$ belongs to $X\cap D$ for all $l\in J_p$, and provided that  for all $j\neq l$,  $\zeta'_1\varphi_p({\zeta'_1}^{\frac1k}\omega^l,\zeta'_2,\ldots,\zeta'_{n-1})\neq \zeta'_1\varphi_p({\zeta'_1}^{\frac1k}\omega^j,\zeta'_2,\ldots,\zeta'_{n-1})$.

For all $\zeta'\in\pi^{-1}_p(\kbz p)$ and all $j\in J_p$, Corollary \ref{cor2.5} and Proposition \ref{prop2.3} give that  $\pi^{-1}_p\big(\zeta'_1,\ldots,\zeta'_{n-1},\zeta'_1\varphi_p({\zeta'_1}^{\frac1k}\omega^l,\zeta'_2,\ldots,\zeta'_{n-1})$ belongs to $X\cap \kb{cC^4|\rho(z)|}z$ which is included in $X\cap D$ provided $c$ is small enough.

Now if we have $\zeta'_1\varphi_p({\zeta'_1}^{\frac1k}\omega^l,\zeta'_2,\ldots,\zeta'_{n-1})= \zeta'_1\varphi_p({\zeta'_1}^{\frac1k}\omega^j,\zeta'_2,\ldots,\zeta'_{n-1})$ for $j\neq l$, then $\pi_p^{-1}\big({\zeta'_1},\zeta'_2,\ldots,\zeta'_{n-1},\zeta'_1\varphi_p({\zeta'_1}^{\frac1k}\omega^l,\zeta'_2,\ldots,\zeta'_{n-1})\big)$ belongs to $\sing(X)$ and so $\zeta'_1=0$. Therefor, it suffices to prove that $H_p$ is bounded when $\zeta'_1\neq 0$ in order to prove that $H_p$ is holomorphic everywhere.

The function $H_p$ is a polynomial in $\zeta'_n$ with coefficients depending on $\zeta'_1,\ldots,\zeta'_{n-1}$ which interpolate $h\circ \pi_p^{-1}$ at the points 
$\big(\zeta'_1,\ldots,\zeta'_{n-1},\zeta'_1\varphi_p({\zeta'_1}^{\frac1k}\omega^l,\zeta'_2,\ldots,\zeta'_{n-1})\big)$, $l\in J_p$. It can therefor be rewritten using divided differences. We set for $\zeta'_1\neq 0$ and $j_1,\ldots, j_l\in J_p$ pairwise distinct:
\begin{align*}
 a_{j_1}(\zeta'_1,\ldots,\zeta'_{n-1})&=
 h\circ \pi_p^{-1}(\zeta'_1,\ldots,\zeta'_{n-1},\zeta'_1\varphi_p({\zeta'_1}^{\frac1k}\omega^{j_1},\zeta'_2,\ldots,\zeta'_{n-1})),\\
 a_{j_1,\ldots,j_{l}}(\zeta'_1,\ldots,\zeta'_{n-1})&=
 \frac{a_{j_1,\ldots,j_{l-1}}(\zeta'_1,\ldots,\zeta'_{n-1})-a_{j_2,\ldots,j_{l}}(\zeta'_1,\ldots,\zeta'_{n-1})}
 {\zeta'_1\varphi_p({\zeta'_1}^{\frac1k}\omega^{j_1},\zeta'_2,\ldots,\zeta'_{n-1})-\zeta'_1\varphi_p({\zeta'_1}^{\frac1k}\omega^{j_{l}},\zeta'_2,\ldots,\zeta'_{n-1})}.
\end{align*}
Writing $J_p$ as $J_p=\{j_1,\ldots,j_{l}\},$ where $l=\# J_p$, we thus have 
\begin{align*}
 &H_p(\zeta'_1,\ldots,\zeta'_n)\\
 &=a_{j_1}(\zeta'_1,\ldots,\zeta'_{n-1})+\ldots+a_{j_1,\ldots,j_{l}}(\zeta'_1,\ldots,\zeta'_{n-1}) \prod_{m=1}^{l-1}
\left(\zeta'_n-\varphi_p({\zeta'_1}^{\frac1k}\omega^{j_{m}},\zeta'_2,\ldots,\zeta'_{n-1})\right).
\end{align*}
We show that for all $j_1,\ldots,j_m\in J_p$ pairwise distinct and all $\zeta'\in\pi_p^{-1}(\kbz p),$ $\zeta'_1\neq 0$, we have $|a_{j_1,\ldots,j_m}(\zeta'_1,\ldots,\zeta'_{n-1})|\leqs |\rho(p)|^{-\frac{m-1}2}$. 
Since for all $j\in J_p$ and $\zeta\in\kbz p$ we have by Proposition \ref{prop2.3} and Corollary \ref{cor2.5}  $\left|\zeta'_n-\varphi_p({\zeta'_1}^{\frac1k}\omega^{j},\zeta'_2,\ldots,\zeta'_{n-1})\right| \leq C^2 |c\rho(p)|^\frac12$, we will have proved that $H_p$ is bounded.

Let $\zeta'$ be a point in $\pi^{-1}_p(\kbz p)$. We consider the disc $\gamma:t\mapsto \pi_p^{-1}(\zeta'_1,\ldots,\zeta'_{n-1},\zeta'_n+t)$ for $t\in D(0,2C^2|c\rho(p)|^{\frac12})$. For all such $t$, $\gamma(t)=\pi_p^{-1}(\zeta')+tw_p$ and since $w_p$ is tangent, $\gamma(t)$ belongs to $\kb{5C^4|c\rho(p)|}p$ which is included in $D$ provided $c$ is small enough. Therefor, by hypothesis, $h\circ\gamma$ admits an holomorphic extension $h_\gamma$ to $D(0,2C^2|c\rho(p)|^{\frac12})$ bounded by some constant which does not depend on $\gamma$.

We set for $j\in J_p$, $t_j=\zeta'_1\varphi_p({\zeta'_1}^{\frac1k}\omega^j,\zeta'_2,\ldots,\zeta'_{n-1})-\zeta'_n$. We then have $|t_j|\leq C^2|c\rho(p)|^{\frac12}$ and
\begin{align*}
 a_{j_1,\ldots,j_m}(\zeta'_1,\ldots,\zeta'_{n-1})&=h_\gamma[t_{j_1},\ldots,t_{j_m}].
\end{align*}
We then get from \cite{Mon} :
\begin{align*}
 a_{j_1,\ldots,j_m}(\zeta'_1,\ldots,\zeta'_{n-1})&=\frac1{2i\pi}\int_{|t|= \frac32C^2|c\rho(p)|^{\frac12}} 
 \frac{h_\gamma(t)}{(t-t_{j_1})\ldots(t-t_{j_m})}dt
\end{align*}
from which we deduce
\begin{align*}
\left| a_{j_1,\ldots,j_m}(\zeta'_1,\ldots,\zeta'_{n-1})\right|&\leqs |\rho(p)|^{-\frac{m-1}2} \sup_{D(0,2cC^2|\rho(p)|^{\frac12})} |h_\gamma|\\
&\leqs |\rho(p)|^{-\frac{m-1}2}
 .\qed
\end{align*}
\subsection{Case $q<+\infty$}
\begin{lemma}\label{lemme4.2}
 Let $\left(\kb{\frac c2|\rho(z_j)|}{z_j}\right)_{j\in\nn}$ be a $\frac c2$-covering of $X\cap D\cap \cu(bD)$ where $\cu(bD)$ is a neighborhood of $bD$. Let $h\in\co(X\cap D)$ be such that for all $j$,  there exists an holomorphic extension $h_j\in\co(\kb{c|\rho(z_j)|}{z_j}\cap L^q (\kb{c|\rho(z_j)|}{z_j})$ of $h$ such that $\sum_{j=0}^{+\infty} \|h_j\|_{L^{q}(\kb{c|\rho(z_j)|}{z_j})}^q$ is finite.\\
 Then there exists an extension $\tilde h\in C^\infty(D)$ of $h$ which satisfies hypothesis (\ref{th0i}-\ref{th0iii}) of Theorem \ref{th0}.
\end{lemma}
\pr 
We define $\tilde h$ as in Lemma \ref{lemme4.1} using our given $\frac c2$-covering and we use the notation of Lemma \ref{lemme4.1}. Thus the Hypothesis (\ref{th0iii}) of Theorem \ref{th0} is immediately satisfied.\\
Cauchy's inequalities then give
\begin{align*}
 \left|
 \diffp{^{|\alpha|}h_j}{{\zeta^*}^\alpha}(\zeta)
 \right|&\leqs \frac1{|\rho(z_j)|^{\alpha_1+\frac{\alpha_2+\ldots+\alpha_n}2 +\frac{n+1}q}} \|h_j\|_{L^{q}(\kb{c|\rho(z_j)|}{z_j})}
\end{align*}
uniformly for all $j$, all multi-indexes $\alpha$ and all $\zeta\in\kb{\frac c2|\rho(z_j)|}{z_j}.$

Therefore, for all multi-indexes $\alpha$ and $\beta$, 
$\left|
\diffp{^{|\alpha|+|\beta|}\psi_jh_j}{{\zeta^*}^\alpha\partial\overline{\zeta^*}^\beta}(\zeta)
\right|\leqs\frac{\|h_j\|_{L^{q}\left(\kb{c|\rho(z_j)|}{z_j}\right)}}{|\rho(z_j)|^{\alpha_1+\beta_1+\frac{\sum_{l=2}^n\alpha_l+\beta_l}2+\frac{n+1}q}}.$
It follows immediately that for all $N\in\nn$, $\rho^{N+n+2} \tilde h$ vanishes to order $N$ on $bD$.

It also comes 
\begin{align*}
\int_{\kb{\frac c2|\rho(z_j)|}{z_j}}\left({|\rho(\zeta)|^{\alpha_1+\frac{\sum_{l=2}^n\alpha_l}2}} \diffp{^{|\alpha|}\psi_jh_j}{\overline{\zeta^*}^\alpha}(\zeta)\right)^{q} d\lambda(\zeta)
&\leqs   \|h_j\|_{L^{q}(\kb{c|\rho(z_j)|}{z_j})}^q.
\end{align*}
Since $\sum_{j=1}^{+\infty} \|h_j\|^q_{L^{q}(\kb{c|\rho(z_j)|}{z_j})}$ is finite, this implies that $\tilde h$ satisfies hypothesis (\ref{th0ii}) of Theorem \ref{th0}.\qed

\begin{theorem}\label{th2}
Let $\left(\kb{\frac c2|\rho(z_j)|}{z_j}\right)$ be a $\frac c2$-covering of $X\cap D\cap \cu(bD)$ where $\cu(bD)$ is a neighborhood of $bD$ and let $h$ belongs to $\co(X\cap D)$ such that for all $j$, there exists $c_j>0$ such that for all holomorphic disc $\gamma:\Delta\to D$ with $\gamma(\Delta)$  included in $\kb{5cC^4|\rho(z_j)|}{z_j}$, there exists $h_\gamma\in\co(\Delta)$ which satisfies
 \begin{enumerate}
  \item $h_\gamma(t)=h\circ\gamma(t)$ for all $t\in\gamma^{-1}(X)$,
  \item $|h_\gamma(t)|^q\leqs\frac1{{\rm Vol}(\kbz{z_j})}c_j,$ for all $t$,
  \item $\sum_{j=0}^{+\infty} c_j$ is finite.
 \end{enumerate}
Then there exists an extension $H$ of $h$ in $\co(D)\cap L^{q}(D)$.
\end{theorem}
\begin{remark}
 As for Theorem \ref{th1}, not all disc have to be tested but only regular disc, i.e.\ disc $\gamma$ such that $\gamma'$ does not vanish. In fact, as we shall see in the proof of Theorem \ref{th2}, the discs that we use are in fact linear discs.

 The conditions of Theorem \ref{th2} are in fact necessary and sufficient. Indeed, if the function $h$ belongs to $\co(D)\cap L^q(D)$, then for any covering $\left(\kb{\frac c2|\rho(z_j)|}{z_j}\right)$ and any disc $\gamma$ included in $\kb{5cC^4|\rho(z_j)|}{z_j}$ for some $j$, we can set $h_\gamma=h\circ \gamma$. Then Cauchy's inequalities will give $|h_\gamma(t)|^q\leqs\frac1{{\rm Vol}(\kbz{z_j})}c_j,$ with $c_j=\|h\|_{L^q(\kb{5cC^4|\rho(z_j)|}{z_j})}^q$. Since any point $z\in D$ belongs to a finite number $M$ of Koranyi balls of a given covering, $M$ which does not depend on $z$, the sum $\sum_{j=0}^{+\infty} \|h\|_{L^q(\kb{5cC^4|\rho(z_j)|}{z_j})}^q$ is bounded up to a multiplicative constant by $\|h\|^q_{L^q(D)}$. 
 
 If $X$ is a manifold, we know from the work of Ohsawa that a function $h\in L^2(X\cap D)\cap \co(X\cap D)$ as a holomorphic extension in $L^2(D)$ and thus satisfies the hypothesis of Theorem \ref{th2}. This can also be checked directly by extending trivially $h$ over a holomorphic disc and by using Cauchy's inequalities in order to prove that this extension satisfies the required estimates.
 \end{remark}

{\it Proof of Theorem \ref{th2} :} In order to apply Lemma \ref{lemme4.2}, we proceed as in the proof of Theorem \ref{th1} from which we use the notations. For $j\in\nn$, we exhibit a holomorphic extension  $h_j\in\co(\kb{c|\rho(z_j)|}{z_j})\cap L^q (\kb{c|\rho(z_j)|}{z_j})$ of $h$ such that$\|h_j\|_{L^{q}(\kb{c|\rho(z_j)|}{z_j})}^q\leqs c_j$, uniformly with respect to $j$.

We put $p=z_j$ and $c(p)=c_j$. 
If $d(p,\sing(X))\leq 10C|c\rho(p)|^{\frac12}$, we set $J_p=\{0,\ldots,k-1\}$ and for any $\xi$, we choose an arbitrary complex number $\xi^{\frac1k}$ such that $\xi=(\xi^{\frac1k})^k$.

If $d(p,\sing(X))\geq 10C|c\rho(p)|^{\frac12}$, we set as in Subsection \ref{Local parametrization} :
$$J_p=\left\{
j\in\{0,\ldots,k\}\tq \exists \zeta' \in\pi_p\big(\kbz p\big),\ \big|\zeta'_1\varphi_p({\zeta'_1}^{\frac1k}\omega^j,\zeta'_2,\ldots,\zeta'_{n-1})-\zeta'_n\big|\leq 2|c\rho(p)|^{\frac12}
\right\}$$
where $\xi\mapsto \xi^{\frac1k}$ is a $k$-root defined on $D(\pi_{p,1}(p),C|c\rho(p)|^{\frac12})$ and $\omega=e^{\frac{2i\pi}k}$.
We then define $h_p$ by setting $h_p=H_p\circ \pi_p^{-1}$ and
\begin{align*}
 H_{p}(\zeta')
 =&\sum_{j\in J_z} \left( \prod_{\over{l\in J_p}{l\neq j}}
\left( \frac{\zeta'_n-\zeta'_1\varphi_p({\zeta'_1}^{\frac1k}\omega^l,\zeta'_2,\ldots,\zeta'_{n-1})}
 {\zeta'_1\varphi_p({\zeta'_1}^{\frac1k}\omega^j,\zeta'_2,\ldots,\zeta'_{n-1})-\zeta'_1\varphi_p({\zeta'_1}^{\frac1k}\omega^l,\zeta'_2,\ldots,\zeta'_{n-1})}\right)\right.\\
 &\hphantom{ \sum_{j\in J_z} \left( \prod_{\over{l\in J_p}{l\neq j}}\right.}\cdot\left. 
 h  \Big(\pi^{-1}_p\big(\zeta'_1,\ldots,\zeta'_{n-1},\zeta'_1\varphi_p({\zeta'_1}^{\frac1k}\omega^l,\zeta'_2,\ldots,\zeta'_{n-1})\big)\Big)
 \vphantom{\prod_{\over{l\in J_p}{l\neq j}}
\left( \frac{\zeta'_n-\zeta'_1\varphi_p({\zeta'_1}^{\frac1k}\omega^l,\zeta'_2,\ldots,\zeta'_{n-1})}
 {\zeta'_1\varphi_p({\zeta'_1}^{\frac1k}\omega^j,\zeta'_2,\ldots,\zeta'_{n-1})-\zeta'_1\varphi_p({\zeta'_1}^{\frac1k}\omega^l,\zeta'_2,\ldots,\zeta'_{n-1})}\right)}\right)
\end{align*}
for $\zeta'\in \pi^{-1}_p(\kb{c|\rho(p)|} p)$.
We rewrite $H_p$ as in the proof of Theorem \ref{th1}. We consider the disc $\gamma:t\mapsto \pi_p^{-1}(\zeta'_1,\ldots,\zeta'_{n-1},\zeta'_n+t)$ for $t\in D(0,2C^2|c\rho(p)|^{\frac12})$. For all such $t$, $\gamma(t)=\pi_p^{-1}(\zeta')+tw_p$ and since $w_p$ is tangent, $\gamma(t)$ belongs to $\kb{5cC^2|\rho(p)|}p$ which is included in $D$ provided $c$ is small enough. Therefor, by hypothesis, $h\circ\gamma$ admits an holomorphic extension $h_\gamma$ to $D(0,2C^2|c\rho(p)|^{\frac12})$ such that
$|h_\gamma|^q\leqs\frac1{{\rm Vol}(\kbz{p})}c(p).$

As in proof of Theorem \ref{th1}, we write $J_p$ as $J_p=\{j_1,\ldots, j_l\}$ and set for $j\in J_p$, $t_j=\zeta'_1\varphi_p({\zeta'_1}^{\frac1k}\omega^j,\zeta'_2,\ldots,\zeta'_{n-1})-\zeta'_n$. We then have $|t_j|\leq C^2|c\rho(p)|^{\frac12}$ and
we get from \cite{Mon} :
\begin{align*}
 a_{j_1,\ldots,j_m}(\zeta'_1,\ldots,\zeta'_{n-1})&=\frac1{2i\pi}\int_{|t|= \frac32C^2|c\rho(p)|^{\frac12}} 
 \frac{h_\gamma(t)}{(t-t_{j_1})\ldots(t-t_{j_m})}dt
\end{align*}
for all $\zeta'\in \pi_p(\kb{c|\rho(p)|} p)$.
It then comes
\begin{align*}
\left| a_{j_1,\ldots,j_m}(\zeta'_1,\ldots,\zeta'_{n-1})\right|&\leqs
\frac1{2\pi|\rho(p)|^{\frac{m-1}2} } \int_0^{2\pi}|h_\gamma(\frac32C^2|c\rho(p)|e^{i\theta})|d\theta.
\end{align*}
Using Jensen's inequality, we get for all $\zeta'\in \pi_p(\kb{c|\rho(p)|} p)$
\begin{align*}
&\left| a_{j_1,\ldots,j_m}(\zeta'_1,\ldots,\zeta'_{n-1})\prod_{j=j_1,\ldots, j_{m-1}}\left(\zeta'_n-\zeta'_1\varphi_p\left({\zeta'_1}^{\frac1k} \omega^j,\zeta'_2,\ldots,\zeta'_{n-1}\right)\right)\right|^q\\
&\leqs
\frac1{2\pi} \int_0^{2\pi}\left|h_\gamma(\frac32cC^2|\rho(p)|e^{i\theta})\right|^qd\theta\\
&\leqs \frac1{{\rm Vol}(\kbz p)}c(p).
\end{align*}
Therefor $H_p$ is bounded in $\{\zeta'\in \pi_p\left(\kbz p\right)/\ \zeta'_1\neq0\}$ and can be extended holomorphically in $\pi_p\left(\kbz p\right)$. Moreover, we have 
\begin{align*}
 \int_{\pi_p\left(\kbz p\right)}|H_p(\zeta')|^qd\lambda(\zeta)&\leqs c(p).\qed
\end{align*}
\begin{remark}
In the proof of Theorems \ref{th1} and \ref{th2}, we have use divided differences in order to construct a holomorphic extension of $h\in\co(X\cap D)$ and we proved that these divided differences are bounded in a certain sense using the existence of good extensions of $h$ in holomorphic discs. But when we look at finding those good extensions in holomorphic discs, a good way (only way ?) is to use divided differences. And in fact, the existence of good extensions in holomorphic discs is equivalent to the control of the divided differences that we used.
\end{remark}

\bibliographystyle{siam}
\bibliography{bibliographie}
\end{document}